\documentclass[12pt]{article}

\usepackage{latexsym}
\usepackage{calc}
\usepackage{datetime}
\usepackage{amssymb, amsmath, amsfonts, listings, mathrsfs} 
\usepackage{xcolor} 
\usepackage[normalem]{ulem} 

\usepackage{graphicx}
\usepackage[labelfont=bf]{caption}

\newcounter{hours}\newcounter{minutes}

\newcommand{\stkout}[1]{\ifmmode\text{\sout{\ensuremath{#1}}}\else\sout{#1}\fi} 

\def\nr{\par \noindent}

\def\Def{\stackrel{\mathrm{def}}{=}}

\def\Ind{\rm Ind \,}

\def\prox{{\rm prox \,}}

\def\dom{{\rm dom \,}}

\def\beq{\begin{equation}}
\def\eeq{\end{equation}}

\def\R{\mathbb{R}}
\def\E{\mathbb{E}}

\def\BI{\begin{itemize}}
\def\EI{\end{itemize}}

\newcommand{\SetEQ}{\setcounter{equation}{0}}
\newcommand{\refLE}[1]{\ensuremath{\stackrel{(\ref{#1})}{\leq}}}
\newcommand{\refEQ}[1]{\ensuremath{\stackrel{(\ref{#1})}{=}}}
\newcommand{\refGE}[1]{\ensuremath{\stackrel{(\ref{#1})}{\geq}}}
\newcommand{\refGT}[1]{\ensuremath{\stackrel{(\ref{#1})}{>}}}

\newcommand{\argmin}[1]{\ensuremath{
\underset{#1}{\arg\min\;}}}

\newtheorem{theorem}{Theorem}
\newtheorem{lemma}{Lemma}
\newtheorem{corollary}{Corollary}

\newtheorem{assumption}{Assumption}
\newtheorem{definition}{Definition}

\newtheorem{example}{Example}
\newtheorem{remark}{Remark}
\newcommand{\proof}{\bf Proof: \rm \nr}
\newcommand{\qed}{\hfill $\Box$ \nr \medskip}
\newcommand{\half}{\mbox{${1 \over 2}$}}

\def\ba{\begin{array}}
\def\ea{\end{array}}
\def\beann{\begin{eqnarray*}}
\def\eeann{\end{eqnarray*}}
\def\bea{\begin{eqnarray}}
\def\eea{\end{eqnarray}}

\def\BT{\begin{theorem}}
\def\ET{\end{theorem}}
\def\BL{\begin{lemma}}
\def\EL{\end{lemma}}
\def\BC{\begin{corollary}}
\def\EC{\end{corollary}}
\def\BE{\begin{example}}
\def\EE{\end{example}}
\def\BD{\begin{definition}}
\def\ED{\end{definition}}
\def\BR{\begin{remark}}
\def\ER{\end{remark}}
\def\BAS{\begin{assumption}}
\def\EAS{\end{assumption}}
\def\BI{\begin{itemize}}
\def\EI{\end{itemize}}

\def\BMP{\begin{minipage}{9.5cm}}
\def\EMP{\end{minipage}}
\def\MPT{\begin{minipage}{11.5cm}}
\def\EPT{\end{minipage}}

\def\la{\langle}
\def\ra{\rangle}

\def\QF{\hspace{5ex} \Box}
\def\QR{\hfill \Box}

\title{High-Order Reduced-Gradient Methods for Composite Variational Inequalities}

\author{Yurii Nesterov
\thanks{Center for Operations Research and Econometrics (CORE),
Catholic University of Louvain (UCL). E-mail:
Yurii.Nesterov@uclouvain.be. This paper has received funding from the European Research Council (ERC) under the European Union’s Horizon 2020 research and innovation program (grant agreement No 788368). The author acknowledges the excellent working conditions at The Hong Kong Polytechnic University, which were important for finalizing the paper. 
}
}

\date{
November 25, 2023
}

\begin{document}
\maketitle

\begin{abstract} 
This paper can be seen as an attempt of rethinking the {\em Extra-Gradient Philosophy} for solving Variational Inequality Problems.
We show that the properly defined {\em Reduced Gradients} can be used instead for finding approximate solutions to Composite Variational Inequalities by the higher-order schemes. Our methods are optimal since their performance is proportional to the lower worst-case complexity bounds for corresponding problem classes. They enjoy the provable hot-start capabilities even being applied to minimization problems. The primal version of our schemes demonstrates a linear rate of convergence under an appropriate uniform monotonicity assumption.
\end{abstract}

\section{Motivation}\label{sc-Motiv}
\SetEQ

Importance of Extra-Gradients in Optimization is commonly justified by their wide employment for solving Variational Inequalities. However, they are also useful for development different
hot-start optimization schemes. For traditional first-order methods, reduction of the distance to the optimum during the minimization process is achieved almost for free. Even the simplest Gradient Method for minimizing convex functions benefits from this property, provided that we choose its steps small enough. However, if the problem is more general (e.g. finding a convex-concave saddle point), then it is not true anymore. 

Similar difficulties experience the basic higher order methods (starting from degree two), where we can only prove that the minimization sequence belongs to a bounded level set of the objective function, defined by the starting point. Note that, these two characteristics, the distance from the starting point to the optimum and the size of the corresponding level set, can be fundamentally different (see Example \ref{ex-LSet}). 
In this paper, we will show that our new Reduced Gradient Framework, which can be seen as amelioration of the concept of Extra-Gradient, can handle these difficulties for different problem settings.

Let us start from historical remarks, related to our
main point of interest, development of efficient methods for solving {\em Variational Inequalities Problems} (VIPs). Below, we use a standard notation, which is anyway explained in the forthcoming Section \ref{sc-Gen}. 

The first fundamental result on Variational Inequalities and Equilibrium Problems is the famous John von Neumann Lemma \cite{JVN}. It justifies the following duality relation:
\beq\label{eq-JVN}
\ba{rcl}
\min\limits_{x \in \Delta_n} \max\limits_{y \in \Delta_m} \la A x, y \ra & =  & \max\limits_{y \in \Delta_m} \min\limits_{x \in \Delta_n}  \la A x, y \ra,
\ea
\eeq
where $\Delta_k = \{ x \in \R^k_+: \sum\limits_{i=1}^k x^{(i)} = 1\}$ and $A \in \R^{m \times n}$. During the next fifty years, this equality served as one of the main motivations for further development of Convex Analysis and Duality Theory, accomplished in the works by Minkowsky, Carath\'eodory, Fenchel, Rockafellar, and others. It forms a basis for the concept of equilibrium in Game Theory and Modern Economics.

Several attempts of highlighting the intrinsic nature of the identity (\ref{eq-JVN}) led to the notion of {\em Variational Inequality} \footnote{An interested reader can find an excellent and very detailed historical survey in the paper \cite{LJ}}. Let $V(\cdot)$ be an operator, defined on a closed convex set $Q \subseteq \E$ with values in $\E^*$. Then Variational Inequality Problem consists in finding a point $x_* \in Q$, which satisfies the following variational principle \cite{Minty}:
\beq\label{def-VarW}
\ba{rcl}
\la V(x), x - x_* \ra & \geq & 0, \quad \forall x \in Q.
\ea
\eeq
Point $x_*$ is called the {\em weak solution} of VIP. Its strong solution $x^* \in Q$ is defined as follows \cite{Stamp}:
\beq\label{def-VarS}
\ba{rcl}
\la V(x^*), x - x^* \ra & \geq & 0, \quad \forall x \in Q.
\ea
\eeq
A standard assumption, ensuring existence of these solutions, is {\em monotonicity} of $V(\cdot)$:
\beq\label{eq-Mon}
\ba{rcl}
\la V(x) - V(y), x - y \ra & \geq & 0, \quad \forall x,y \in Q.
\ea
\eeq
In this case, if operator $V(\cdot)$ is continuous, then $x_* = x^*$. 

Note that for $\E = \R^n$, definition (\ref{def-VarS}) can be written in a fixed-point form:
\beq\label{def-VarP}
\ba{rcl}
x^* & = & \pi_Q(x^* - \lambda V(x^*)),
\ea
\eeq
where $\lambda > 0$ is an arbitrary fixed parameter, and $\pi_Q(y)$ denotes a Euclidean projection of point $y$ onto the set $Q$ in the standard norm $\| \cdot \|_{(2)}$.

Let us mention two important sources of Variational Inequalities.
\BI
\item{\underline{\sc Optimization problems}}. Consider the problem 
\beq\label{prob-Min}
f^* \; = \; \min\limits_{x \in Q} f(x), 
\eeq
where $f(\cdot)$ is a convex and continuously differentiable function on $Q$. Then its optimal solution $x^* \in Q$ can be characterised by the variational principle (\ref{def-VarS}) with $V(x) = \nabla f(x)$. In this case, inequality (\ref{def-VarS}) is just the well known first-order optimality condition.
\item{\underline{\sc Saddle-point problems}}.
Let $Q_1 \subseteq \E_1$ and $Q_2 \subseteq \E_2$ be two closed and convex sets. Assume that for any $x \in Q_1$ function $\phi(x,\cdot)$ is continuously differentiable and concave, and for any $y \in Q_2$, function $\phi(\cdot,y)$ is continuously differentiable and convex. Then, under mild assumptions (e.g. Sion's Theorem \cite{Sion}), we have
\beq\label{def-Sad}
\ba{rcl}
\min\limits_{x \in Q_1} \max\limits_{y \in Q_2} \phi(x,y) & = & \max\limits_{y \in Q_2} \min\limits_{x \in Q_1}\phi(x,y).
\ea
\eeq
The solution $(x^*,y^*) \in Q_1 \times Q_2$ of this problem is characterized by the following condition:
$$
\ba{rcl}
\phi(x^*,y) & \leq & \phi(x^*,y^*) \; \leq \phi(x,y^*), \quad \forall x \in Q_1,\; y \in Q_2.
\ea
$$
At the same time, this pair is also a solution of the Variational Inequality (\ref{def-VarS}) with $Q = Q_1 \times Q_2$ and $V(x,y) = (\nabla_1 \phi(x,y), - \nabla_2 \phi(x,y)) \in \E^*_1 \times \E^*_2$.

For the {\em Bilinear Matrix Game} (\ref{eq-JVN}), we have $Q = \Delta_n \times \Delta_m$ and 
$$
\ba{rcl}
V(x,y) & = & (A^Ty, -Ax) \Def C z \in \R^{n+m}, 
\ea
$$
where $z = (x,y) \in \R^{n+m}$ and $C \in \R^{(n+m) \times (n+m)}$.
Note that matrix $C$ is {\em skew-symmetric}, and hence {\em positive-semidefinite}:
\beq\label{def-PSD}
\ba{rcl}
\la C z, z \ra & \equiv & 0, \quad \forall z \in \R^{n+m}.
\ea
\eeq
\EI
 
However, from the computational point of view, an attractive generality of  formulation~(\ref{def-VarS}) is sometimes misleading. Indeed, it appears that the class of optimization problems is {\em provably easier} than the general saddle-point problems (and consequently, than the general VIPs). 

In order to see the origin of the difficulties, let us look at behaviour of the simplest Gradient Method in both situations. Firstly, consider the problem of convex minimization with simple set constraint:
\beq\label{prob-Un}
\ba{rcl}
f(x^*) & = & \min\limits_{x \in Q} f(x), 
\ea
\eeq
where $Q \subseteq \R^n$ is a closed convex set, and function $f(\cdot)$ has Lipschitz continuous gradients:
\beq\label{def-LipG}
\ba{rcl}
\| \nabla f(x) - \nabla f(y) \|_{(2)} & \leq & L \| x - y \|_{(2)}, \quad x, y 
\in \R^n.
\ea
\eeq
Let us choose in the Gradient Method 
\beq\label{met-GM}
\ba{rcl}
x_{k+1} & = & \pi_Q(x_k - h_k \nabla f(x_k)), \quad k \geq 0, 
\ea
\eeq
the stepsizes $h_k \in (0,{1 \over L})$. Note that the projection $x_{k+1}$ is characterized by the following variational principle:
\beq\label{eq-VGM}
\ba{rcl}
\la x_{k+1} - x_k + h_k \nabla f(x_k), x - x_{k+1} \ra & \geq & 0, \quad \forall x \in Q.
\ea
\eeq
Therefore, denoting $r_k = \| x_k - x^* \|_{(2)}$, we get
$$
\ba{rcl}
r_{k+1}^2 & = &  r_k^2 + 2 \la x_{k+1} - x_k, x_{k+1} - x^* \ra - \| x_{k+1}-x_k \|^2_{(2)}\\ \\
& \refLE{eq-VGM} & r_k^2 + 2 h_k \la \nabla f(x_k),x^* - x_k \ra - 2h_k \left[ \la \nabla f(x_k),x_{k+1}- x_k \ra + {1 \over 2h_k} \| x_{k+1}-x_k \|^2_{(2)} \right]\\ \\
& \refLE{def-LipG} & r_k^2 + 2 h_k [f(x^*) - f(x_k)] - 2 h_k [ f(x_{k+1}) - f(x_k)] \; = \; r_k^2 - 2 h_k [f(x_{k+1}) - f(x^*)] .
\ea
$$
Thus, in method (\ref{met-GM}), the distance to the optimal point is {\em monotonically decreasing}. Moreover, the constant step strategy
\beq\label{eq-CStep}
\ba{rcl}
h_k & = & {1 \over L}, \quad k \geq 0,
\ea
\eeq
results in the following rate of convergence:
\beq\label{eq-GMRate}
\ba{rcl}
f(x_k) - f(x^*) & \leq & {L r_0^2\over 2 k}, \quad k \geq 1.
\ea
\eeq

Consider now a variant of the Gradient Method, as applied to Variational Inequality (\ref{def-VarS}):
\beq\label{met-GVI}
\ba{rcl}
x_{k+1} & = & \pi_Q( x_k - h_k V(x_k)), \quad k \geq 0,
\ea
\eeq
with some $h_k > 0$. Then,
\beq\label{eq-ProgVI}
\ba{rcl}
r_{k+1}^2 & \leq & r_k^2  - 2 h_k \la V(x_k) , x_k - x^* \ra + h_k^2 \| V(x_k) \|_{(2)}^2.
\ea
\eeq
However, now we can have $\la V(x_k) , x_k - x^* \ra = 0$ (see (\ref{def-PSD})). Thus, for example, for skew-symmetric operators, like $V(x) = Cx$ with $C+C^T=0$ and $x^* = 0$, the distance to the solution can only {\em increase}. Hence, the sequence $\{ x_k \}_{k \geq 0}$, in general, {\em cannot} converge to the solution.

A big difference in the complexity of Optimization Problems and VIs can be observed also in non-degenerate case. Let us assume that the objective function in (\ref{prob-Un}) is strongly convex with parameter $\mu > 0$. Then, applying method (\ref{met-GM}) with $h_k = {2 \over \mu + L}$, we have (see, for example, Theorem 2.2.14 in \cite{LN}):
\beq\label{eq-RateUn}
\ba{rcl}
r_{k+1}^2 & \leq \left({L-\mu \over L+\mu} \right)^2 r_k^2.
\ea
\eeq

For VI with Lipschitz-continuous and strongly monotone operator, 
\beq\label{eq-VIClass}
\ba{rcl}
\| V(x) - V(y) \|_{(2)} & \leq & L \| x - y \|_{(2)}, \\
\\
\la V(x) - V(y) , x - y \ra & \geq & \mu \| x - y \|_{(2)}^2, \quad \forall x, y \in \R^n,
\ea
\eeq
we get the following bound:
$$
\ba{rcl}
r_{k+1}^2 & \stackrel{(\ref{def-VarP}),(\ref{met-GVI})}{\leq} & \| x_{k} - x^* - h_k(V(x_k) - V(x^*)) \|_{(2)}^2\\
\\
& = & r_k^2 - 2 h_k \la V(x_k) - V(x^*), x_k - x^* \ra + h_k^2 \| V(x_k) - V(x^* )\|_{(2)}^2\\
\\
& \refLE{eq-VIClass} & \left(1 - 2 \mu h_k + h_k^2 L^2 \right) r_k^2.
\ea
$$
Thus, choosing the optimal stepsize $h_k = {\mu \over L^2}$, we get the following rate of convergence:
\beq\label{eq-RateVIS}
\ba{rcl}
r_{k+1}^2 & \leq \left(1 - \left({\mu \over L}\right)^2 \right) r_k^2.
\ea
\eeq
Clearly, it is much worse than the rate (\ref{eq-RateUn}).

From the discussion above, we can see that, as compared to optimization problems, the simplest Gradient Method encounter significant difficulties in controlling the distance to the solution of Variational Inequality. 

Recently, it was observed that the similar difficulties experience the basic higher order methods. 
Consider, for example, the Cubic Newton Method \cite{CN}
$$
\ba{rcl}
x_{k+1} & = & \arg\min\limits_{x \in \R^n} \Big\{ \la \nabla f(x_k), x - x_k \ra + \half \la \nabla^2 f(x_k)(x-x_k), x - x_k \ra + {M \over 6} \| x - x_k \|^3_{(2)} \Big\},
\ea
$$
as applied to problem (\ref{prob-Min}) with $Q = \R^n$, where the objective function has Lipschitz continuous Hessian:
\beq\label{def-LipH}
\ba{rcl}
\| \nabla^2 f(x) - \nabla^2 f(y) \| & \leq & M \| x - y \|_{(2)}, \quad x, y \in \R^n.
\ea
\eeq
In accordance to Theorem 4.2.2 in \cite{LN}, this method has the following rate of convergence:
\beq\label{eq-QNRate}
\ba{rcl}
f(x_k) - f(x^*) & \leq & {9 M D_0^3 \over (k+4)^2}, \quad k \geq 1,
\ea
\eeq
where $D_0 = \max\limits_x \{ \| x - x^* \|_{(2)}: \; f(x) \leq f(x_0)\}$. However, this diameter can be very big even if the initial distance to the optimum is small. Let us look at the following example.
\BE\label{ex-LSet}
Consider the following function:
$$
\ba{rcl}
f(x) & = & |x^{(1)}|^3 + \sum\limits_{i=1}^{n-1}|x^{(i+1)} - 2 x^{(i)}|^3, \quad x \in \R^n,
\ea
$$
with $x^* = 0$ and $f(x^*) = 0$.
Then for the point $x_0 = (1, \dots, 1)^T$, we have $f(x_0) = n$. Similarly, for the point $\bar x \in \R^n$ with coordinates
$$
\ba{rcl}
\bar x^{(i)} & = & 2^i - 1, \quad i = 1, \dots , n,
\ea
$$
we have $f(\bar x) = n$. At the same time, $\| x_0 \|^2_{(2)} = n$ and 
$$
\ba{rcl}
\| \bar x \|^2_{(2)} & = & {8 \over 3}(2^n-1)(2^{n-1}-1)+n. 
\ea
$$
Thus, for this function, the constant $D_0$ in the right-hand side of inequality (\ref{eq-QNRate}) is exponentially big. \qed
\EE

The main goal of this paper consists in developing a regularization technique, which helps to control the distance between the generated test points and the solution. At each iteration of the corresponding schemes, this technique needs computation of an additional gradient (or, an operator value for VIs). It results in the optimal high-order methods for solving Variational Inequalities. Moreover, for corresponding methods, it ensures the {\em hot-start} capabilities.

\section{Reduced-Gradient Cut for Vatiational Inequalities}\label{sc-Extra}
\SetEQ

Before discussing the idea of Reduced Gradient, let us look at the theoretical guarantees for the Gradient Method (\ref{met-GVI}), where the monotone operator $V(\cdot)$ is Lipschitz continuous on the convex set $Q$. However, in order to avoid unlimited increase of the distance to the solution, we assume now that the set $Q$ is bounded:
\beq\label{ass-QBound}
\ba{rcl}
\| x - x^* \|_{(2)} & \leq & D \; < \; + \infty, \quad \forall x \in Q.
\ea
\eeq

As before, the point $x_{k+1}$ in definition (\ref{met-GVI}) satisfies the following variational principle:
\beq\label{eq-VVI}
\ba{rcl}
\la x_{k+1} - x_k + h_k V(x_k), x - x_{k+1} \ra & \geq & 0, \quad \forall x \in Q.
\ea
\eeq
Denote $r_k(x) = \| x - x_k \|_{(2)}$ and $\rho_k = \| x_k - x_{k+1} \|_{(2)}$. Then,
\beq\label{eq-RecVI}
\ba{rcl}
r_{k+1}^2(x) & = &  r_k^2(x) + 2 \la x_{k+1} - x_k, x_{k+1} - x \ra - \rho_k^2\\ \\
& \refLE{eq-VVI} & r_k^2(x) + 2 h_k \la V(x_k),x - x_{k+1} \ra -  \rho_k^2\\
\\
& \refLE{eq-VIClass} & r_k^2(x) + 2 h_k \la V(x_{k+1}),x - x_{k+1} \ra + 2 h_k L D \rho_k -  \rho_k^2\\
\\
& \refLE{eq-Mon} & r_k^2(x) + 2 h_k \la V(x), x - x_{k+1} \ra +  h^2_k L^2 D^2.
\ea
\eeq

The last inequality suggests the following {\em merit function} for a candidate solution $\bar x \in Q$:
\beq\label{def-MF}
\ba{rcl}
\mu_Q(\bar x) & = & \sup\limits_{x \in Q} \; \la V(x), \bar x - x \ra.
\ea
\eeq
Note that function $\mu_Q(\cdot)$ is convex and non-negative on $Q$. Moreover, $\mu_Q(\bar x) = 0$ if and only if the point $\bar x$ is a weak solution of VI (see (\ref{def-VarW})). Since we assume $D < +\infty$, this merit function is well defined for all $\bar x \in Q$. However, if the operator $V(\cdot)$ has a potential, then its merit function is well define on $Q$ even without this additional assumption. 
\BL\label{lm-FRes}
Let $V(x) = \nabla f(x)$ for some convex function $f(\cdot)$. Then, for any $\bar x \in Q$, we have
\beq\label{eq-FRes}
\ba{rcl}
\mu_Q(\bar x) & \leq & f(\bar x) - f^*.
\ea
\eeq
\EL
\proof
Indeed, $\mu_Q(\bar x) = \sup\limits_{x \in Q} \; \la \nabla f(x), \bar x - x \ra \leq \sup\limits_{x \in Q} \; [f(\bar x) - f(x)] = f(\bar x) - f^*$.
\qed

Thus, in the potential case, $\mu_Q(\cdot)$ is a weaker optimality measure than the usual functional residual.

Now, let us use the relation (\ref{eq-RecVI}) in order to justify the quality of some {\em average points} of the sequence $\{ x_k \}_{k \geq 0}$, generated by method (\ref{met-GVI}). For the segment of iterations $[p,m]$, denote
$$
\ba{rcl}
S^1_{p,m} & = & \sum\limits_{i=p}^{m} h_i, \quad S^2_{p,m} \; = \; \sum\limits_{i=p}^{m} h_i^2, \quad x_{p,m} \; = \; {1 \over S^1_{p,m}}\sum\limits_{i=p}^{m} h_i x_{i+1}.
\ea
$$

\BL\label{lm-Seg}
For any $0 \leq p \leq m$, we have
\beq\label{eq-Seg}
\ba{rcl}
\mu_Q(x_{p,m}) & \leq & {1 + L^2 S^2_{p,m} \over 2 S^1_{p,m}} D^2.
\ea
\eeq
\EL
\proof
Indeed, summing up inequalities (\ref{eq-RecVI}) for $k = p, \dots, m$, we get
$$
\ba{rcl}
2 S^1_{p,m} \la V(x), x_{p,m} - x \ra + r_{m+1}^2(x) & \leq & r_p^2(x) + L^2 D^2 S^2_{p,m} \; \leq \; (1 + L^2 S^2_{p,m})D^2,
\ea
$$
and (\ref{eq-Seg}) follows.
\qed

The analytical structure of the right-hand side of inequality (\ref{eq-Seg}) often arise in analysis of subgradient methods for Nonsmooth Optimization (e.g. Section 3.2.3 in \cite{LN}). For approaching its minimal value, it is recommended to choose
\beq\label{eq-GStep}
\ba{rcl}
h_k & = & {1 \over L \sqrt{k+1}}, \quad k \geq 0.
\ea
\eeq
In this case, we have
$$
\ba{rcl}
L^2 S^2_{0,m} & = & 1 + \sum\limits_{k=1}^m {1 \over k+1} \; \leq \; 1 + \int\limits_0^m {dt \over t+1} \; = \; 1 + \ln (m+1),\\
\\
S^1_{0,m} & = & {1 \over L}\sum\limits_{k=0}^m {1 \over \sqrt{k+1}} \; \geq \; {1 \over L} \int\limits_0^{m+1} {dt \over \sqrt{t+1}} \; = \; {2 \over L}(\sqrt{m+2}-1).
\ea
$$
Hence, using the upper bound (\ref{eq-Seg}), we get the following rate of convergence:
\beq\label{eq-GMLog}
\ba{rcl}
\mu_Q(x_{0,m}) & \leq & {2 + \ln (m+1) \over 4(\sqrt{m+2}-1)} L D^2.
\ea
\eeq

An unpleasant logarithmic factor in this bound can be eliminated by {\em sliding averaging} (see Section 4 in \cite{Arik}). Indeed, for any $m \geq 1$, we have
$$
\ba{rcl}
L^2 S^2_{m,2m-1} & = & \sum\limits_{k=m}^{2m-1} {1 \over k+1} \; \leq \; \int\limits_{m-1}^{2m-1} {dt \over t+1} \; = \; \ln 2,\\
\\
S^1_{m,2m-1} & = & {1 \over L}\sum\limits_{k=m}^{2m-1} {1 \over \sqrt{k+1}} \; \geq \; {1 \over L} \int\limits_m^{2m} {dt \over \sqrt{t+1}} \; = \; {2 \over L}(\sqrt{2m+1}-\sqrt{m+1}) \; \geq \; {2 \sqrt{m-1} \over (1+\sqrt{2})L} .
\ea
$$
Thus, for any $m \geq 2$, we obtain
\beq\label{eq-GMVRate}
\ba{rcl}
\mu_Q(x_{m,2m-1}) & \leq & {(1+\sqrt{2})\ln 2 \over 2 \sqrt{m-1}} L D^2.
\ea
\eeq

Despite to convergence results (\ref{eq-GMLog}), (\ref{eq-GMVRate}), we have to acknowledge that the trajectory of the Gradient Method (\ref{met-GVI}) for VIPs can be very strange. Indeed, let $V(x) = Cx$ for a skew-symmetric matrix $C$, and $\psi(\cdot) = \mbox{Ind}_Q(\cdot)$ for a bounded closed convex set $Q$, containing zero as an interior point. Then $x^* = 0$ and the distance to the solution along the trajectory $\{ x_k \}_{k \geq 0}$ is going to increase up to the limits prescribed by the boundary of the set $Q$ (see \ref{eq-ProgVI})). Thus, it will go around the origin along an orbit with potentially big diameter. In this situation, only the average points can converge to the solution.

At the same time, as we have seen, the rates (\ref{eq-GMLog}), (\ref{eq-GMVRate}) are much worse as compared to the rate (\ref{eq-GMRate}) of the Gradient Method being applied to minimization problem (\ref{prob-Un}). Moreover, the latter rate is obtained for a stronger optimality measure (see Lemma \ref{lm-FRes}). 

The above comparison explains the importance of the following {\em Extra-Gradient Method}:
\beq\label{met-ExGrad}
\ba{rcl}
y_k & = & \pi_Q(x_k - h V(x_k)), \quad x_{k+1} \; = \; \pi_Q(x_k - h V(y_k)), \quad k \geq 0,
\ea
\eeq
with some $h > 0$, which was proposed in \cite{Korp} (see also \cite{Arik} for the modern justification, interpretations, and applications). It appears that for Lipschitz continuous monotone operators the method (\ref{met-ExGrad}) ensures the following convergence rate:
$$
\ba{rcl}
\mu_Q(\tilde x_k) & \leq & O\left({1 \over k}\right), \quad k \geq 1,
\ea
$$
where $\tilde x_k$ are some average points.
Moreover, this rate is unimprovable for the corresponding problem class (see Section 5 in \cite{Arik}).

In this paper, we present an alternative to the extragradient idea, the framework of {\em Reduced Gradients}. It works well both for optimization problems and variational inequalities. It can be applied for justifying methods of any order.

Our framework is based on two important elements.
The first one is the sequence of prox-centers $\{ v_t \}_{t \geq 0}$, with a controlled distance to the optimum. And the second one is the main sequence $\{x_t \}_{t \geq 0}$, which justifies the global rate of convergence of the process and helps in updating the prox-centers in a regular way.

Let us introduce this framework for the following {\em Composite Variational Inequality} (CVI):
\beq\label{prob-CVar}
\ba{c}
\mbox{Find $x^* \in \dom \psi: \; \la V(x), x - x^* \ra + \psi(x) \geq \psi(x^*), \; \forall x \in \dom \psi$},
\ea
\eeq
where $\psi(\cdot)$ is a proper closed convex function and $V(\cdot): \dom \psi \to \E^*$ is a monotone operator.

An important variant of this problem corresponds to a 
variational inequality with monotone operator $V(\cdot)$, defined on a closed convex set $Q$. Then we can chose $\psi(\cdot) = \Ind_Q(\cdot)$, the indicator function of $Q$. Another variant is the {\em composite minimization problem} \cite{Comp}, with $V(\cdot)$ being a gradient of a differentiable convex function defined on a closed convex set $\dom \psi$. Then inequality (\ref{prob-CVar}) is just a necessary and sufficient characterization of the optimal point $x^* \in \dom \psi$.

Let us introduce an equivalent formulation of CVI.
\BL\label{lm-CVar1}
Let operator $V(\cdot)$ be continuous and monotone on $\dom \psi$.
The point $x^* \in \dom \psi$ is a solution to problem (\ref{prob-CVar}) if and only if 
\beq\label{prob-CVar1}
\ba{rcl}
\la V(x^*), x - x^* \ra + \psi(x) & \geq & \psi(x^*), \quad \forall x \in \dom \psi.
\ea
\eeq
\EL
\proof
Indeed, in view of monotonicity of $V(\cdot)$, inequality (\ref{prob-CVar1}) implies (\ref{prob-CVar}). Let us assume now that (\ref{prob-CVar}) is valid for all $x \in \dom \psi$. Denote $x_{\alpha} = x^* + \alpha(x-x^*)$, $0 \geq \alpha \leq 1$. Then,
$$
\ba{rcl}
\psi(x^*) & \refLE{prob-CVar} & \la V(x_{\alpha}) , x_{\alpha} - x^* \ra + \psi(x_{\alpha}) \; \leq \; \alpha \la V(x_{\alpha}), x - x^* \ra + (1 - \alpha) \psi(x^*) + \alpha \psi(x).
\ea
$$
Dividing this inequality by $\alpha > 0$, we get
$$
\ba{rcl}
0 & \leq & \la V(x_{\alpha}), x - x^*\ra +  \psi(x) - \psi(x^*).
\ea
$$
It remains to take the limit as $\alpha \to 0$.
\qed

Following the tradition, we call the point $x^*$ satisfying condition (\ref{prob-CVar}) the {\em weak solution} of the corresponding CVI, and (\ref{prob-CVar1}) defines its {\em strong solution}. The strong CVI (\ref{prob-CVar1}) can be posed in the following form:
\beq\label{prob-Inc}
\mbox{\fbox{$\;$ Find $x^* \in \dom \psi: \; -V(x^*) \in \partial \psi(x^*)\;$}}
\eeq

In the first part of this paper (Sections \ref{sc-Primal}, \ref{sc-DExtra}, and \ref{sc-Proj}), we analyse the iterative schemes in the most general form, assuming only a possibility of performing an \fbox{\em Essential Step} from some point $v \in \dom \psi$. This implies an ability of finding a point $T = T(v) \in \dom \psi$ such that
\beq\label{def-Ess}
\ba{rcl}
\la V_{\psi}(T), v - T \ra & > & 0, 
\ea
\eeq
for $V_{\psi}(T) = V(T) + \psi'_T$ with some $\psi'_T \in \partial \psi(T)$. Latter on, we will discuss different ways for computing such points and describing their quality. 

An important characteristic of the Essential Step is the \fbox{\em Universal Stepsize} defined as 
\beq\label{def-ASize}
\ba{rcl}
a(v,T) & \Def & {1 \over \| V_{\psi}(T) \|_*^2} \la V_{\psi}(T), v - T \ra.
\ea
\eeq
This value is used for performing the Reduced-Gradient steps in {\em all methods} below, independently on the problem classes containing the operator $V(\cdot)$.

More specifically, all our strategies are as follows. Given the current prox-center $v_t$, we compute the next point $x_{t+1} = T(v_t)$ and the {\em reduced gradient} $V_{\psi}(x_{t+1}$), which {\em cut} the prox-center from the solution set $X^*$:
\beq\label{eq-Cut}
\ba{rcl}
\la V_{\psi}(x_{t+1}), v_t - x_{t+1} \ra & \refGT{def-Ess} & 0 \; \refGE{prob-CVar} \; \la V_{\psi}(x_{t+1}), x^* - x_{t+1} \ra , \quad \forall x^* \in X^* \subseteq \dom \psi.
\ea
\eeq
And then, we generate a new prox-center $v_{t+1}$ with a smaller distance to the localization set
\beq\label{def-QT}
\ba{rcl}
Q_t & \Def & \Big\{ x \in \dom \psi: \; \la V_{\psi}(x_{t+1}), x_{t+1} - x \ra \geq 0 \Big\}.
\ea
\eeq
Since $X^* \subseteq Q_t$, we manage to control the distance
from the points $\{ v_t \}_{t \geq 0}$ to the solution set.

In order to make our methods efficient, we need to generate {\em deep cuts}, which ensure a big difference between the left- and right-hand sides of inequality (\ref{eq-Cut}). Note that for optimization methods we can often take $x_{t+1} = v_t$ since in the case of potential operator $V(x) = \nabla f(x)$ the deepness of the cut can be easily ensured by different standard inequalities like
$$
\ba{rcl}
\la \nabla f(x_{t+1}), x_{t+1} - x^* \ra & \geq & f(x_{t+1}) - f(x^*), \quad \la \nabla f(x_{t+1}), x_{t+1} - x^* \ra \; \geq \; {1 \over L} \| \nabla f(x) \|_*^2, \\
\\
\la \nabla f(x_{t+1}), x_{t+1} - x^* \ra & \geq & \mu \| x_{t+1} - x^* \|^2, \quad \mbox{etc.}
\ea
$$

The situation with Variational Inequalities is much more delicate. Indeed, as we have seen, it is possible to have $\la V(x_{t+1}), x_{t+1} - x^* \ra = 0$. Hence, our main hope consists in developing some technique, which ensures a big universal stepsize $a(v,T)>0$. Later on, we will show how this can be done for different problem classes. At this point, we just mention two general results, which we often use for guaranteeing the validity of the second part of inequality (\ref{eq-Cut}).

Let operator $A(\cdot)$ be monotone on $\dom \psi$. Denote by $T = T_{\psi}(A) \in \dom \psi$ 
the strong solution of the following CVI:
\beq\label{eq-CondA}
\ba{rcl}
\la A(T), x - T \ra + \psi(x) & \geq & \psi(T), \quad \forall x \in \dom \psi.
\ea
\eeq
\BL\label{lm-Cut2}
For $T = T_{\psi}(A)$, {\em Reduced Gradient} \fbox{$V_{\psi}(T) \Def V(T) - A(T)$} ensures
\beq\label{eq-RGLow}
\ba{rcl}
\la V_{\psi}(T), T - x \ra & \geq & \la V(T), T - x \ra + \psi(T) - \psi(x), \quad x \in \dom \psi.
\ea
\eeq
Hence, for $x = x^*$, we get a correct cut:
$\la V_{\psi}(T), x^* - T \ra \refLE{prob-CVar} 0$.
\EL
\proof
Indeed, for any $x \in \dom \psi$, we have
$$
\ba{rcl}
\la V_{\psi}(T), T - x \ra & = & \la V(T) - A(T), T - x \ra \; \refGE{eq-CondA} \; \la V(T), T - x \ra + \psi(T) - \psi(x). \QF
\ea
$$

\BL\label{lm-VPsi}
For monotone operators $A_i(\cdot)$ denote $T_i = T_{\psi}(A_i)$, $i = 1, 2$. Then
\beq\label{eq-VPsi}
\ba{rcl}
\la V_{\psi}(T_1) - V_{\psi}(T_2), T_1 - T_2 \ra & \geq & \la V(T_1) - V(T_2), T_1 - T_2 \ra \; \geq \; 0
\ea
\eeq
\EL
\proof
Indeed, adding two versions of inequality (\ref{eq-RGLow}) with $T=T_1$, $x = T_2$, and $T=T_2$, $x = T_1$, we get (\ref{eq-VPsi}).
\qed

The remaining part of this paper is organized as follows. In Section \ref{sc-Gen}, we introduce notations and prove some general auxiliary results. Sections \ref{sc-Primal}, \ref{sc-DExtra}, and \ref{sc-Proj} are devoted to descriptions of three general reduced-gradent methods (primal, dual, and projective ones), for which we can prove general convergence results in terms of the size of the universal steps. In the next three Sections \ref{sc-Rates}, \ref{sc-Var}, and \ref{sc-Uni} we prove the lower bounds for the universal step sizes, which justify the rate of convergence of the corresponding schemes. These are the high-order reduced-gradient methods for minimization problems, high-order reduced-gradient methods for VIPs, and high-order reduced-gradient methods for uniformly convex VIPs. In all cases, we prove the rates of convergence for the corresponding merit functions (e.g. residual in the function value) and/or for the measures of optimality condition (e.g. norm of the reduced gradients). We finish the paper with Conclusion, were we summarize our rezultes and discuss the related papers \cite{ABJS,LJ}.

\section{Notations and generalities}\label{sc-Gen}
\SetEQ

In this paper, we consider problems with variables belonging to a finite-dimensional real vector space $\E$. Its dual space $\E^*$ is formed by all linear functions on $\E$. The value of $s \in \E^*$ at point $x \in \E$ is denoted by $\la s, x \ra$. 

For a convex function $f(\cdot)$, we denote by $\partial f(x)$ its subdifferential at point $x \in \dom f$.
For differentiable function $f(\cdot)$ with open domain $\dom f$, we denote by $\nabla f(x) \in \E^*$ its gradient at point $x \in \dom f$:
$$
\ba{rcl}
f(y) & \approx & f(x) + \la \nabla f(x), y - x \ra, \quad y \in \dom f.
\ea
$$
Notation $\nabla^2 f(x)$ is used for the Hessian of $f(\cdot)$ at $x$. Note that the Hessian is a self-adjoint linear operator from $\E$ to $\E^*$:
$$
\ba{rcl}
\nabla f(y) & \approx & \nabla f(x) + \nabla^2f(x)(y-x) \; \in \; \E^*, \quad y \in \dom f.
\ea
$$

Let us consider nonlinear operators, which map an open subset of $\E$ into another space $\E_0$. In this paper, we are mainly interested in two cases, $\E_0 = \R$, which corresponds to minimization problems, and $E_0 = E^*$, which corresponds to VIPs.

For an operator $V(\cdot): \dom V \to \E_0$, denote by
$$
\ba{c}
D^pV(x)[h_1, \dots, h_p] \; \in \E_0, \quad p \geq 0,
\ea
$$
the $p$th derivative of operator $V(\cdot)$ at point $x \in \dom V$ along directions $h_i \in \E$, $1 \leq i \leq p$. 
Note that the value $D^pV(x)[h_1, \dots, h_p]$ is not changing after arbitrary permutation of directions. In the case $h_i = h$, $1 \leq i \leq p$, we use a shorter notation $D^pV(x)[h]^p$. It is convenient to treat the value  $D^0V(x)[h]^0$ as a constant operator $V(x)$.

Nonlinear operators with $\E_0 = \R$ are just real-valued functions. In this case,
$$
\ba{rcl}
D^0f(x)[h]^0 & = & f(x), \quad
Df(x)[h] \; = \; \la \nabla f(x), h \ra, \quad D^2f(x)[h]^2 \; = \; \la \nabla^2 f(x)h, h \ra.
\ea
$$

For the important family of operators with $\E_0 = \E^*$, the first derivative $DV(x)$ (called {\em Jacobian}) belongs to the class of linear operators from $\E$ to $\E^*$ (notation ${\cal L}(\E)$):
$$
\ba{rcl}
V(y) & \approx & V(x) + DV(x)[y-x], \quad y \in \dom V.
\ea
$$
This feature allows us to introduce an important class of {\em monotone operators}, for which
\beq\label{def-Mon}
\ba{rcl}
\la DV(x)h, h \ra & \geq 0, \quad x \in \dom V, \; h \in \E,
\ea
\eeq
(notation $DV(x) \succeq 0$).
If $V(x) = \nabla f(x)$ for some function $f(\cdot)$, then $DV(x) = \nabla^2 f(x)$, and condition (\ref{def-Mon}) implies convexity of function $f(\cdot)$.

The directional derivatives are used for defining the {\em Taylor series}, which can predict the value of the map $V(\cdot)$ at other points. Denote
$$
\ba{rcl}
{\cal T}_{x,p}^V(y) & = & \sum\limits_{i=0}^p {1 \over i!} D^iV(x)[y-x]^i \; \in \; \E_0, \quad y \in \E.
\ea
$$
If no ambiguity arise, we omit the superscript $V$ in this notation.

We can expect that ${\cal T}_{x,p}^V(y) \approx V(y)$ when $y$ is close to $x$. Indeed, by Taylor Theorem, assuming continuity of $(p+1)$th derivative, we have
\beq\label{eq-TV}
\ba{rcl}
V(y) & = & {\cal T}^V_{x,p}(y) + \int\limits_0^1 {(1-\tau)^p \over p!} D^{p+1}V(x+\tau(y-x))[y-x]^{p+1} d \tau, \quad y \in \dom V.
\ea
\eeq

For predicting Jacobians, we can use the derivative of Taylor series, that is 
$$
\ba{rcl}
D{\cal T}_{x,p}^V(y) & = & \sum\limits_{i=1}^p {1 \over (i-1)!} D^iV(x)[y-x]^{i-1} \; \Def \;
{\cal J}_{x,p}^V(y) \; \in \; {\cal L}(\E), \quad y \in \E.
\ea
$$
Then again, by Taylor Theorem we have
\beq\label{eq-TD}
\ba{rcl}
DV(y) & = & {\cal J}^V_{x,p}(y) + \int\limits_0^1 {(1-\tau)^{p-1} \over (p-1)!} D^{p+1}V(x+\tau(y-x))[y-x]^{p} d \tau, \quad y \in \dom V.
\ea
\eeq
However, for qualifying the sharpness of approximation (\ref{eq-TV}), (\ref{eq-TD}), we need to employ norms.

For measuring distances in $\E$, we can use a general norm $\| \cdot \|$. Then the dual norm is defined in the standard way,
\beq\label{def-DNorm}
\ba{rcl}
\| g \|_* & = & \max\limits_{x \in \E} \{ \la g, x \ra: \; \| x \| \leq 1 \}, \; g \in \E^*,
\ea
\eeq
ensuring the important Cauchy-Schwartz inequality:
\beq\label{eq-CS}
\ba{rcl}
\la s, x \ra & \leq & \| s \|_*  \cdot \| x \|, \quad x \in \E, \; s \in \E^*.
\ea
\eeq

Let us measure distances in the image space $\E_0$ by some norm $\| \cdot \|_0$. Then we can measure the size of a symmetric polylinear form $A[\cdot, \dots , \cdot ]:\; \E^p \to \E_0$ by the following quantity:
$$
\ba{rcl}
\| A \|_0 & = & \max\limits_{h \in \E} \Big\{ \| A[h]^p \|_0 : \; \| h \| \leq 1 \Big\}.
\ea
$$
If we assume the uniform boundedness of  $(p+1)$st derivative on a convex domain ${\cal D} \subseteq \E$:
\beq\label{eq-DBound}
\ba{rcl}
\| D^{p+1} V(x) \|_0 & \leq & M_{p+1}({\cal D}), \quad \forall x \in {\cal D},
\ea
\eeq
then, it is easy to see that
\beq\label{eq-VDiv}
\ba{rcl}
\| V(y) - {\cal T}^V_{x,p}(y) \| & \refLE{eq-TV} & {M_{p+1}({\cal D}) \over (p+1)!} \| y - x \|^{p+1}, \quad x, y \in {\cal D}.
\ea
\eeq

For measuring changes in Jacobians, we need two more characteristics. Firstly, for $B \in {\cal L}(\E)$, we introduce the following induced norm:
\beq\label{def-JNorms}
\ba{rcl}
\| B \|_{\cal L} & = & \max\limits_{x, y \in \E} \{ \la B x , y \ra: \; \| x \| \leq 1, \; \| y \| \leq 1\} \; \refEQ{def-DNorm} \; \max\limits_{x \in \E} \{ \| B x \|_*: \; \| x \| \leq 1 \}.
\ea
\eeq
And secondly, for a polylinear form $A: \E^p \to \E_0$, we need a stronger norm
$$
\ba{rcl}
\| A \|^0 & = & \max\limits_{h_1, \dots, h_p \in \E} \Big\{ \| A[h_1, \dots, h_p] \|_0: \; \| h_i \| \leq 1, 1 \leq i \leq p \Big\}.
\ea
$$
In general, $\| A \|^0 \geq \| A \|_0$. Now, assuming uniform boundedness of  $(p+1)$st derivative on ${\cal D}$:
\beq\label{eq-D2Bound}
\ba{rcl}
\| D^{p+1} V(x) \|^0 & \leq & \widehat M_{p+1}({\cal D}), \quad \forall x \in {\cal D},
\ea
\eeq
we can see that
\beq\label{eq-V2Div}
\ba{rcl}
\| D V(y) - {\cal J}^V_{x,p}(y) \|_{\cal L} & \refLE{eq-TD} & {\widehat M_{p+1}({\cal D}) \over p!} \| y - x \|^{p}, \quad x, y \in {\cal D}.
\ea
\eeq

In order to use a general norm in our methods, we need to introduce a differentiable convex prox-function $d(\cdot)$, $\dom d \subseteq \E$, which is strongly convex with respect to this norm with parameter one:
\beq\label{def-SConv}
\ba{rcl}
d(y) & \geq & d(x) + \la \nabla d(x), y - x \ra + \half \| y - x \|^2, \quad x, y \in \dom d.
\ea
\eeq
By this function, we can define the {\em Bregman distance} between two points in $\dom d$:
\beq\label{def-Breg}
\ba{rcl}
\beta_d(x,y) & = & d(y) - d(x) - \la \nabla d(x), y - x \ra \;  \refGE{def-SConv} \; \half \| y - x \|^2, \quad x, y \in \dom d.
\ea
\eeq
If no ambiguity arise, we use a shorter notation $\beta(\cdot,\cdot) \equiv \beta_d(\cdot,\cdot)$. Then, the {\em proximal gradient operator} at point $\bar x \in \dom \psi$ is defined as follows:
\beq\label{def-Prox}
\ba{c}
\prox_{\bar x,h} (g) = \argmin{x \in \dom \psi} \Big\{ h \la g, x - \bar x \ra  + \beta(\bar x, x ) \Big\}, \quad g \in \E^*,
\ea
\eeq
where $h \geq 0$ is a {\em step-size parameter}. We assume that functions $d(\cdot)$ and $\psi(\cdot)$ are simple enough for computing the exact value of this operator in a closed form. For $T = \prox_{\bar x, h}(g)$, the first-order optimality condition is as follows:
\beq\label{eq-CondProx}
\ba{rcl}
\la h g + \nabla d(T) - \nabla d(\bar x), x - T \ra & \geq & 0, \quad \forall x \in \dom \psi.
\ea
\eeq

Sometimes, we need to use Euclidean norms. For that, we choose a positive definite self-adjoint linear operator $B = B^*: \E \to \E^*$, and define the norm as follows:
\beq\label{def-Euclid}
\ba{rcl}
\| x \|_B & = & \la B x, x \ra^{1/2}, \quad x \in \E, \quad  \quad \| g \|^*_B \; = \; \la g, B^{-1}g \ra^{1/2}, \quad g \in \E^*.
\ea
\eeq
In this case, we always choose $d(x) = \half \| x \|^2_B$, which results in the following simple form for the Bregman distance:
\beq\label{def-EBreg}
\ba{rcl}
\beta_d(x,y) & = & \half \| x - y \|^2_B , \quad x, y \in \E.
\ea
\eeq
Hence, the proximal gradient operator looks as follows:
\beq\label{def-EProx}
\ba{c}
\prox_{\bar x,h}(g) = \argmin{x \in \dom \psi} \; \half \| x - \bar x + h B^{-1} g \|^2_B.
\ea
\eeq
If $\psi(\cdot)$ is an indicator function of a closed convex set $Q \subseteq \E$, that is
$$
\ba{rcl}
\psi(x) & = & \Ind_Q(x) \; \Def\; \left\{ \ba{rl} 0, & \mbox{if $x \in Q$},\\ + \infty, & \mbox{otherwise,} \ea \right.
\ea
$$
then $\prox_{\bar x,h}(g) = \pi_Q(\bar x - h B^{-1} g)$,
where $\pi_Q(x)$ is a Euclidean projection of point $x \in \E$ onto the set $Q$.

In the higher-order methods, we often need to use functions $d_p(x) = {1 \over p} \| x \|^p_B$ with $p \geq 2$ (see Section~4.2.2 in \cite{LN}). These functions have the following differential characteristics:
\beq\label{eq-DPChar}
\ba{c}
\nabla d_p(x) = \| x \|^{p-2}_B Bx, \quad \nabla^2 d_p(x) = \| x \|^{p-2}_B \Big[ B + (p-2) {Bx x^* B \over \| x \|^2_B} \Big] \succeq \| x \|^{p-2}_B B, \quad x \in \E.
\ea
\eeq
It can be proved that these functions are uniformly convex of degree $p$ \cite{LN}:
\beq\label{eq-DUni}
\ba{rcl}
d_p(y) - d_p(x) - \la \nabla d_p(x), y - x \ra & \geq & {1 \over p} \left( \half \right)^{p-2} \| y - x \|^p_B,\\
\\
\la \nabla d_p(y) - \nabla d_p(x), y - x \ra & \geq & \left( \half \right)^{p-2} \| y - x \|^p_B, \quad x, y \in \E.
\ea
\eeq

Finally, we need the following technical lemma.
\BL\label{lm-Tech}
For any $\sigma \geq 1$,  $\gamma >0$, and $\delta \geq 0$, we have
\beq\label{eq-Tech}
\ba{rcl}
\inf\limits_{g > 0} \Big\{ \half \gamma g^{2 \over \sigma} + g^{1 - \sigma \over \sigma} \delta \Big \} & = & {\sigma + 1 \over 2} \left( {\gamma \over \sigma - 1} \right)^{\sigma - 1 \over \sigma + 1} \delta^{2 \over 1+\sigma}.
\ea
\eeq
\EL
\proof
Indeed, if $\sigma = 1$, this equality is evident. Consider the case $\sigma > 1$. Note that after the change of variables $\tau = g^{2 \over \sigma}$, the objective function in this minimization problem is convex in $\tau$. Therefore, the optimal $g_* > 0$ can be found from the equation
$$
\ba{rcl}
{\gamma \over \sigma} g^{{2 \over \sigma} - 1}_* & = & {\sigma-1 \over \sigma} g^{{1 \over \sigma} - 2}_* \delta \quad \Leftrightarrow \quad
g_* \; = \; \left[ {\sigma - 1 \over \gamma} \delta \right]^{\sigma \over 1 + \sigma}.
\ea
$$
Hence, the optimal value of the objective in (\ref{eq-Tech}) can be computed as follows:
$$
\ba{rcl}
g^{2 \over \sigma}_* \left( \half \gamma + g^{-{1+\sigma \over \sigma}}_* \delta \right)  & = & g^{2 \over \sigma}_* \left( \half \gamma +  \left[ {\sigma - 1 \over \gamma} \delta \right]^{-1} \delta \right) \;= \; g^{2 \over \sigma}_* \left( \half \gamma +  {\gamma \over \sigma - 1} \right) \\
\\
& = & { \gamma (\sigma + 1) \over 2(\sigma - 1)}  g^{2 \over \sigma}_*. \QR
\ea
$$

\section{Primal Reduced-Gradient Method}\label{sc-Primal}
\SetEQ

The simplest implementation of the ideas discussed in Section \ref{sc-Extra} looks as follows.
\beq\label{met-EPrimal}
\ba{|l|}
\hline \\
\hspace{15ex} \mbox{\bf Primal Reduced-Gradient Method}\\
\\
\hline \\
\ba{rl}
\, & \mbox{{\bf 0.} Choose $x_0 \in \dom \psi$. Set $v_0 = x_0$.}\\
& \\
& \mbox{{\bf 1. $t$-th iteration ($t \geq 0$).}}\\
& \\
& \mbox{{\bf a).} Compute $x_{t+1} \in\dom \psi$ satisfying condition (\ref{eq-Cut}).} \quad \\
& \\
& \mbox{{\bf b).} Set $a_{t+1} \refEQ{def-ASize} a(v_t,x_{t+1})$ and $v_{t+1} = \prox_{v_t,a_{t+1}}(V_{\psi}(x_{t+1}))$.}\\
\ea \\
\\
\hline
\ea
\eeq

The rate of convergence of this scheme, as well as that of all methods below, can be characterized by the following objects:
\beq\label{def-Obj}
\ba{rcl}
b_t & = & {\la V_{\psi}(x_t), v_{t-1} - x_t \ra^2 \over 2 \| V_{\psi}(x_t) \|^2_*} \; = \; \half a^2_t \| V_{\psi}(x_t) \|^2_*,\\
\\
A_0 & = & 0, \quad A_t \; = \; \sum\limits_{i=1}^{t} a_i, \quad
B_0 \; = \; 0, \quad B_t \; = \; \sum\limits_{i=1}^{t} b_i, \\

\\
V_t(x) & = & \la V_{\psi}(x_t), x_t - x\ra, \quad t \geq 1, \quad x \in \E.
\ea
\eeq
For the geometric interpretation of these schemes, note that the value ${1 \over \| V_{\psi}(x_t) \|_*}\la V_{\psi}(x_t), v_{t-1} - x_t \ra$ is equal to the distance from the point $v_{t-1}$ to the hyperplane $\{ x \in \E: \; \la V_{\psi}(x_t),x_t- x  \ra = 0\}$, measured in the norm $\| \cdot \|$.

\BT\label{th-EPrimal}
For any $t \geq 0$ and $x \in \dom \psi$, we have
\beq\label{eq-EPrimal}
\ba{rcl}
\sum\limits_{i=1}^t a_i V_i(x) + B_t + \beta(v_t,x) & \leq & \beta(x_0,x).
\ea
\eeq
\ET
\proof
For $i \geq 0$, in view of inequality (\ref{eq-CondProx}), we have
$$
\ba{rcl}
\la a_{i+1} V_{\psi}(x_{i+1}) + \nabla d(v_{i+1}) - \nabla d(v_i), x - v_{i+1} \ra & \geq & 0
\ea
$$
for all $x \in \dom \psi$.
Therefore,
$$
\ba{rcl}
\beta(v_{i+1},x) & = & d(x) - d(v_{i+1}) - \la \nabla d(v_{i+1}), x - v_{i+1} \ra\\
\\
& \leq & d(x) - d(v_{i+1}) + \la a_{i+1} V_{\psi}(x_{i+1}) - \nabla d(v_i), x - v_{i+1} \ra\\
\\
& = & \beta(v_i,x) + d(v_i) + \la \nabla d(v_i), v_{i+1} - v_i \ra - d(v_{i+1}) \\
\\
& & + a_{i+1}  \la  V_{\psi}(x_{i+1}), x - v_{i+1} \ra \\
\\
& \refLE{def-SConv} & \beta(v_i,x) - \half \| v_{i+1} - v_i \|^2 + a_{t+1} \la  V_{\psi}(x_{i+1}), x - v_{i+1} \ra .
\ea
$$
Thus, maximizing the right-hand side of this inequality in $\| v_i - v_{i+1} \|$, we get
$$
\ba{rcl}
\beta(v_{i+1},x) & \leq & \beta(v_i,x) + \half a_{i+1}^2 \|V_{\psi}(x_{i+1}) \|^2_* + a_{i+1} \la  V_{\psi}(x_{i+1}), x - v_{i} \ra \\
\\
& = & \beta(v_i,x) + \half a_{i+1}^2 \| V_{\psi}(x_{i+1}) \|^2_*  + a_{i+1} [ \la  V_{\psi}(x_{i+1}), x_{i+1} - v_{i} \ra - V_{i+1}(x)]\\
\\
& = & \beta(v_i,x) - b_{i+1}
- a_{i+1}  V_{i+1}(x).
\ea
$$
Summing up this inequalities for $i=0, \dots, t-1$, we get the bound (\ref{eq-EPrimal}) 
\qed

Note that the rules of method (\ref{met-EPrimal}) do not depend on the problem class containing the operator $V(\cdot)$.

In the case of $\psi(\cdot)$ being an indicator function of a closed convex set $Q$, we often chose the norm $\| \cdot \|$ to be Euclidean. Then, the scheme of method (\ref{met-EPrimal}) becomes very simple.

\beq\label{met-EEPrim}
\ba{|l|}
\hline \\
\hspace{3ex} \mbox{\bf Euclidean Primal Reduced-Gradient Method}\\
\\
\hline \\
\ba{rl}
\, & \mbox{{\bf 0.} Choose $x_0 \in \dom \psi$. Set $v_0 = x_0$.}\\
& \\
& \mbox{{\bf 1. $t$-th iteration ($t \geq 0$).}}\\
& \\
& \mbox{{\bf a).} Compute $x_{t+1} \in\dom \psi$ satisfying condition (\ref{eq-Cut}).} \quad \\
& \\
& \mbox{{\bf b).} Set $a_{t+1} = {1 \over (\| V_{\psi}(x_{t+1}) \|^*_B)^2} \la V_{\psi}(x_{t+1}), v_t - x_{t+1} \ra$ and}\\
\\
& \mbox{define $v_{t+1} = \pi_Q\Big(v_t - a_{t+1} B^{-1}V_{\psi}(x_{t+1})\Big)$.}\\
\ea \\
\\
\hline
\ea
\eeq

In Sections \ref{sc-Rates} and \ref{sc-Var}, we will see how to transform inequality (\ref{eq-EPrimal}) into the rates of convergence of method (\ref{met-EPrimal}) for different problem settings. Before that, in the next two sections, we justify two more variants of the reduced-gradient methods.

\section{Dual Reduced Gradient Method}\label{sc-DExtra}
\SetEQ

Dual optimization schemes are usually supported by the machinery of estimating sequences. Consider the following method.

\beq\label{met-DEGM}
\ba{|l|}
\hline \\
\hspace{14ex} \mbox{\bf Dual Reduced Gradient Method}\\
\\ \hline \\
\ba{rl}
& \mbox{{\bf 0).} Choose $x_0 \in \dom \psi$ and define $\Psi_0(x) = \beta(x_0,x)$.}\\
& \\ 
& \mbox{\bf 1). $t$-th iteration ($t \geq 0$).}\\
& \\
& \mbox{{\bf a)} Compute $v_t = \arg\min\limits_{x \in \dom \psi} \Psi_t(x)$.}\\
& \\
& \mbox{{\bf b)} Compute $x_{t+1} \in \dom \psi$ satisfying condition (\ref{eq-Cut}).}\\
& \\
& \mbox{{\bf c)} Choose $a_{t+1} =  a(v_t,x_{t+1})$ and define}\\
& \\
& \quad \Psi_{t+1}(x) = \Psi_t(x)+ a_{t+1} [\la V(x_{t+1}), x - x_{t+1} \ra + \psi(x)].\\
\ea\\
\\
\hline
\ea
\eeq

As compared with method (\ref{met-EPrimal}), the main information is kept now in the dual space. Indeed, if we denote
$s_t= \sum\limits_{i=1}^t a_i V(x_i)$ for $t \geq 0$,
then 
$$
\ba{rcl}
v_t & = & \min\limits_{x \in \dom \psi} \Big[ \la s_t, x \ra + \beta(x_0,x) + A_t \psi(x) \Big],
\ea
$$
Therefore, this point has the following characterization:
\beq\label{eq-DProx}
\ba{rcl}
\la s_t + \nabla d(v_t) - \nabla d(x_0), x - v_t \ra + A_t \psi(x) & \geq & A_t \psi(v_t), \quad \forall x \in \dom \psi.
\ea
\eeq
It is important that $s_t$ aggregates the linear functionals $V(x_i)$, and not $V_{\psi}(x_i)$.

Let us justify the convergence of this method denoting 
$\Psi^*_t = \min\limits_{x \in \dom \psi} \Psi_t(x)$.
\BT\label{th-DEM}
For any $t \geq 0$, we have
\beq\label{eq-DEM}
\ba{rcl}
\sum\limits_{i=1}^t a_i \psi(x_i) + B_t & \leq & \Psi^*_t.
\ea
\eeq
\ET
\proof
For $t =0$, we have $B_0 = 0$ and $\Psi^*_0 = 0$. Hence, the inequality (\ref{eq-DEM}) is valid. Note that
$$
\ba{rcl}
\Psi^*_{t+1} & = & \Psi_t(v_{t+1}) + a_{t+1} [ \la V(x_{t+1}), v_{t+1} - x_{t+1} \ra + \psi(v_{t+1})]\\
\\
& \geq & \Psi^*_t + \half \| v_{t+1} - v_t \|^2 + a_{t+1} [ \la V_{\psi}(x_{t+1}), v_{t+1} - x_{t+1} \ra + \psi(x_{t+1})]\\
\\ 
& \geq & \Psi^*_t - \half a^2_{t+1} \|V_{\psi}(x_{t+1}) \|^2_* + a_{t+1} [ \la V_{\psi}(x_{t+1}), v_t - x_{t+1} \ra + \psi(x_{t+1})]\\
\\ 
& \refEQ{def-Obj} & \Psi^*_t + b_{t+1} + a_{t+1} \psi(x_{t+1}). \QR
\ea
$$

In Sections \ref{sc-Rates} and \ref{sc-Var}, we will show how to transform the inequality (\ref{eq-DEM}) into the rate of convergence of the corresponding schemes.

\section{Projecting Reduced Gradient Method}\label{sc-Proj}
\SetEQ

In this method, the idea of decreasing the distance to the points of the set $Q_t$ (see (\ref{def-QT})) is implemented in the most straightforward way. 

\beq\label{met-PREGM}
\ba{|l|}
\hline \\
\hspace{11ex} \mbox{\bf Projecting Reduced Gradient Method}\\
\\ \hline \\
\ba{rl}
& \mbox{{\bf 0).} Choose $x_0 \in \dom \psi$. Set $v_0 = x_0$.}\\
& \\ 
& \mbox{\bf 1). $t$-th iteration ($t \geq 0$).}\\
& \\
& \mbox{{\bf a)} Compute $x_{t+1} \in \dom \psi$ satisfying condition (\ref{eq-Cut}).}\\
& \\
& \mbox{{\bf b)} Set $v_{t+1}=\argmin{x \in Q_t}\; \beta(v_t,x)$ with $Q_t$ defined by (\ref{def-QT}).}\\
\ea\\
\\
\hline
\ea
\eeq

\BT\label{th-PREGM}
Let the sequences $\{ v_t \}_{t \geq 0}$ and $\{ x_t \}_{t \geq 0}$ be generated by method (\ref{met-PREGM}). Then, for any $t \geq 0$, we have
\beq\label{eq-PREGM}
\ba{rcl}
\beta(v_t,x^*) + B_t & \leq & \beta(x_0,x^*).
\ea
\eeq
\ET
\proof
Indeed, from the first-order optimality condition, defining the point $v_{t+1}$, we have
$$
\ba{rcl}
\la \nabla d(v_{t+1}) - \nabla d(v_t), x - v_{t+1} \ra & \geq & 0
\ea
$$
for all $x \in \dom \psi$ with $\la V_{\psi}(x_{t+1}), x_{t+1} - x \ra \geq 0$. 
Therefore,
$$
\ba{rcl}
\beta(v_{t+1},x^*) & = & d(x^*) - d(v_{t+1}) - \la \nabla d(v_{t+1}), x^* - v_{t+1} \ra \\
\\
& \leq & d(x^*) - d(v_{t+1}) - \la \nabla d(v_{t}), x^* - v_{t+1} \ra\\
\\
& = & \beta(v_t,x^*) + d(v_t) - d(v_{t+1}) - \la \nabla d(v_{t}), v_t - v_{t+1} \ra\\
\\
& = & \beta(v_t,x^*) - \beta(v_t,v_{t+1}) \; \refLE{def-SConv} \; \beta(v_t,x^*) \; - \half \| v_t - v_{t+1} \|^2.
\ea
$$
At the same time, in view of condition (\ref{eq-Cut}), $v_t \not\in Q_t$. Thus,
$$
\ba{rcl}
\la V_{\psi}(x_{t+1}), x_{t+1} - v_{t+1} \ra & = & 0.
\ea
$$
Hence,
$$
\ba{rcl}
b_{t+1} & = & {\la V_{\psi}(x_{t+1}), v_t - x_{t+1} \ra^2 \over 2 \| V_{\psi}(x_{t+1}) \|^2_*} \; = \; {\la V_{\psi}(x_{t+1}), v_t - v_{t+1} \ra^2 \over 2 \| V_{\psi}(x_{t+1}) \|^2_*} \; \leq \; \half \| v_t - v_{t+1} \|^2. \QF
\ea
$$

\section{Rates of convergence for Composite Minimization}\label{sc-Rates}
\SetEQ

Let us show how to transform the statements of Theorems \ref{th-EPrimal}, \ref{th-DEM}, and \ref{th-PREGM} into the convergence rates of the corresponding optimization schemes.
Consider the following problem of Composite Minimization:
\beq\label{prob-COpt}
\min\limits_{x \in \dom \psi} \Big\{ F(x) = f(x) + \psi(x) \Big\},
\eeq
where function $f(\cdot)$ is $p$ times continuously differentiable on $\dom \psi$ and its $p$th derivative is Lipschitz continuous on $\dom \psi$ with constant $L_p$. From now on, we assume that all distances are measured in the Euclidean norm $\| \cdot \| \equiv \| \cdot \|_B$. Hence, $\beta(x_0,x) = \half \| x - x_0\|^2$.

Consider the Taylor polynomial of degree $p \geq 1$ for function $f(\cdot)$:
$$
\ba{rcl}
{\cal T}_{x,p}^f(y) & = & \sum\limits_{i=0}^p {1 \over i!} D^i f(x)[y-x]^i.
\ea
$$
As it was explained in Section \ref{sc-Gen}, for all $x, y \in \dom \psi$, we have
\beq\label{eq-FDiff}
\ba{rcl}
\Big| f(y) - {\cal T}_{x,p}^f(y) \Big| & \refLE{eq-VDiv} & {L_p \over (p+1)!} \| y - x \|^{p+1}.
\ea
\eeq
Applying (\ref{eq-VDiv}) to $V(x) = \nabla f(x)$, we get
\beq\label{eq-GDiff}
\ba{rcl}
\Big\| \nabla f(y) - \nabla {\cal T}_{x,p}^f(y) \Big\|_* & \leq & {L_p \over p!} \| y - x \|^{p}. 
\ea
\eeq

Let us define now the Basic Tensor Step at point $v \in \dom \psi$ as follows:
\beq\label{def-BTStep}
\ba{rcl}
x_+ \; = \; T_M^p(v) & = & \argmin{y \in \dom \psi} \Big\{ {\cal T}_{v,p}^f(y) + {M \over (p+1)!} \| y - v \|^{p+1} + \psi(y) \Big\},
\ea
\eeq
with $M \geq p L_p$. From \cite{Imp} we know that the optimization problem in (\ref{def-BTStep}) is convex. Hence, the operator 
$$
\ba{rcl}
A_{v,p}^f(y) & \Def &  \nabla {\cal T}_{v,p}^f(y) + {M \over p!} \| y - v \|^{p-1}B(y-v)
\ea
$$
is monotone. Therefore, the point
$T = T_M^p(v)$, is a strong solution to the following CVI:
$$
\ba{rcl}
\la A_{v,p}^f(T), y - T \ra + \psi(y) & \geq & \psi(T), \quad y \in \dom \psi.
\ea
$$
Thus, in accordance to Lemma \ref{lm-Cut2}, we can choose
$$
\ba{rcl}
V_{\psi}(T) & = & \nabla f(T) - A_{v,p}^f(T) \; \in \; \partial F(T),
\ea
$$
and establish some lower bounds for the objects (\ref{def-Obj}). Let us prove first the following lemma.

\BL\label{lm-VIStep}
For any $v \in \dom \psi$ and $T = T_M^p(v)$ with $M \geq pL_p$ and $p \geq 1$, we have
\beq\label{eq-VIStep}
\ba{rcl}
\la V_{\psi}(T), v - T \ra  & \geq & \gamma_p (M) \| V_{\psi}(T) \|_*^{p+1 \over p},\\
\\
\gamma_p(M) & \Def & {p \over M} \left( p! \over p+1 \right)^{1 \over p} \left[ M^2 - L_p^2 \over p^2 - 1 \right]^{p-1 \over 2p}.
\ea
\eeq
\EL
\proof
Denote $r = \| T - v \|$. Then
$$
\ba{rcl}
{L_p r^p \over p! } & \refGE{eq-GDiff} & \| \nabla f(T) - \nabla {\cal T}_{v,p}^f(T) \|_* \; = \; \| V_{\psi}(T) + {M r^{p-1} \over p!} B(T - v) \|_*.
\ea
$$
Squaring both sides of this inequality and rearranging the terms, we get
$$
\ba{rcl}
2 {M r^{p-1} \over p!} \la V_{\psi}(T), v - T \ra & \geq & \| V_{\psi}(T) \|_*^2 + {M^2 - L_p^2 \over (p!)^2} r^{2p}.
\ea
$$
Hence,
$2 \la V_{\psi}(T), v - T \ra \geq {p! \over M r^{p-1}} \| V_{\psi}(T) \|_*^2 + {M^2 - L_p^2 \over (p!) M} r^{p+1}$. Let us use the inequality
$$
\ba{rcl}
a r^{1-p} + b r^{1+p} & \geq & 2 p \left[ a \over p+1 \right]^{p+1 \over 2p} \left[ b \over p-1 \right]^{p-1 \over 2p}, \quad a, b \geq 0,
\ea
$$ 
which can be proved by minimizing its left-hand side in $r > 0$. Taking
$$
\ba{rcl}
a & = & {p! \over M} \| V_{\psi}(T) \|_*^2, \quad b \; = \; {M^2 - L_p^2 \over p! M},
\ea
$$
we get
$$
\ba{rcl}
2 \la V_{\psi}(T), v - T \ra & \geq & 2 p \left[ {p! \over (p+1) M} \| V_{\psi}(T) \|_*^2 \right]^{p+1 \over 2p} \left[ {M^2 - L_p^2 \over (p-1) p! M} \right]^{p-1 \over 2p} \\
\\
& = & {2p \over M} \left[ p! \over p+1 \right]^{1 \over p} \left[ M^2 - L_p^2 \over p^2 - 1 \right]^{p-1 \over 2p} \| V_{\psi}(T) \|_*^{p+1 \over p}. \QF
\ea
$$

Inequality (\ref{eq-VIStep}) ensures that
$$
\ba{rcl}
a_t & \refGE{def-ASize} & \gamma_p(M) \| V_{\psi}(T) \|_*^{1 - p \over p},\quad
b_t \; \refGE{def-Obj} \; \half \gamma_p^2(M) \| V_{\psi}(T) \|_*^{2 \over p}.
\ea
$$
Let us look what does it mean for method (\ref{met-EPrimal}) as applied to problem (\ref{prob-COpt}) with the Essential Step Rule (\ref{def-BTStep}). In view of inequality (\ref{eq-EPrimal}), for any $x \in \dom \psi$, we have
\beq\label{eq-InterP}
\ba{rcl}
\half \| x_0 - x \|^2 \; \geq \; B_t + \sum\limits_{i=1}^t a_i \la V_{\psi}(x_i), x_i - x \ra & \geq & B_t + \sum\limits_{i=1}^t a_i [F(x_i) - F(x)] \\
& = & B_t + A_t [\tilde F_t - F(x)],
\ea
\eeq
where $\tilde F_t = {1 \over A_t} \sum\limits_{i=1}^t a_i F(x_i)$.
Taking $x = x^*$ and denoting $g_t = \| V_{\psi}(x_i)\|_*$, we get
\beq\label{eq-FTilde}
\ba{rcl}
\half \| x_0 - x^* \|^2 & \geq & \sum\limits_{i=1}^t \Big\{ \half \gamma_p^2(M) (g_i)^{2 \over p} + \gamma_p(M) g_i^{1-p \over p} [\tilde F_t - F_*] \Big\}.
\ea
\eeq 

Since $\tilde F_t - F_* \geq 0$, for $g_t^* = \min\limits_{1 \leq i \leq t} g_i$, we have
\beq\label{eq-RateGP}
\ba{rcl}
g_t^* & \leq & \left[ {1 \over \gamma_p(M)} \| x_0 - x^* \| \right]^p \left(1 \over t \right)^{p \over 2}, \quad t \geq 1.
\ea
\eeq
Note that the function $\gamma_p(M)$ is decreasing in $M \geq \sqrt{p} L_p$. Hence, its minimal value is achieved at $M^*_p = p L_p$. That is 
$$
\ba{rcl}
\gamma^*_p & = & \gamma_p(M^*_p) \; = \; \left( p! \over (p+1) L_p \right)^{1 \over p}.
\ea
$$
Thus, we get the following guarantee for method (\ref{eq-EPrimal}) with $M = M^*_p$:
\beq\label{eq-RateGP*}
\ba{rcl}
g_t^* & \leq &   {(p+1) L_p \over p!} \| x_0 - x^* \|^p \left(1 \over t \right)^{p \over 2}, \quad t \geq 1.
\ea
\eeq

Further, applying to each term in the right-hand side of inequality (\ref{eq-FTilde}) the lower bound (\ref{eq-Tech}) with $\gamma = \gamma_p(M)$, $\sigma = p$, $g=g_i$, and $\delta = \tilde F_t - F_*$, we get
$$
\ba{rcl}
\half \| x_0 - x^* \|^2 & \geq & t \cdot \gamma_p(M) {p+1 \over 2} \left( \gamma_p(M) \over p-1 \right)^{p-1 \over p+1} [\tilde F_t - F_*]^{2 \over p+1}.
\ea 
$$ 
This means that
\beq\label{eq-RateFP}
\ba{rcl}
\tilde F_t - F_* & \leq & {\| x_0 - x^* \|^{p+1} \over p+1} \left({1 \over \gamma_p(M)}\right)^p \left(p-1 \over p+1\right)^{p-1 \over 2} \left(1 \over t \right)^{p+1 \over 2}.
\ea
\eeq
Finally, using $M = M^*_p$, we come to the following bound:
\beq\label{eq-RateFP*}
\ba{rcl}
\tilde F_t - F_* & \leq & {L_p \| x_0 - x^* \|^{p+1} \over p!} 
\left(p-1 \over p+1 \right)^{p-1 \over 2}  \left(1 \over t \right)^{p+1 \over 2}.
\ea
\eeq

The main advantage of the estimates (\ref{eq-RateGP*}) and (\ref{eq-RateFP*}) consists in their \underline{\em hot-start capability}: the closeness of the starting point to the optimum simplifies computation of an approximate solution in function value.

Relation (\ref{eq-InterP}) can be used for justifying an accuracy certificate of method (\ref{met-EPrimal}). For that, we need to know an upper bound for the distance from the initial point $x_0$ to the set of optimal solutions:
$$
\ba{rcl}
\| x_0 - x^* \| & \leq & R_0
\ea
$$
for some $x^* \in X^*$. Let us define the following value:
\beq\label{def-ACerfF}
\ba{rcl}
\Delta^F_t & = & {1 \over A_t}\max\limits_{x \in \dom \psi} \Big\{ \sum\limits_{i=1}^t a_i \la V_{\psi}(x_i), x_i - x \ra: \; \| x - x_0 \| \leq R_0 \Big\} \\
\\
& \geq & {1 \over A_t}\max\limits_{x \in \dom \psi} \Big\{ \sum\limits_{i=1}^t a_i [F(x_i) - F(x)]: \; \| x - x_0 \| \leq R_0 \Big\} \; = \; \tilde F_t - F_*.
\ea
\eeq
At the same time, $\half R_0^2 \refGE{eq-InterP} B_t + A_t \Delta^F_t$. Hence, using the same arguments as for bounding the gap $\tilde F_t - F_*$, we come to the following guarantees:
\beq\label{eq-RateCF}
\ba{rcl}
\Delta^F_t & \leq & {R_0^{p+1} \over p+1} \left({1 \over \gamma_p(M)}\right)^p \left(p-1 \over p+1\right)^{p-1 \over 2} \left(1 \over t \right)^{p+1 \over 2},
\ea
\eeq
which is valid for all $M \geq p L_p$. 
For $M = M^*_p$, it can be specified as follows:
\beq\label{eq-RateCF*}
\ba{rcl}
\Delta^F_t & \leq & {L_p R_0^{p+1} \over p!} \left({p-1 \over p+1} \right)^{p-1 \over 2}
 \left(1 \over t \right)^{p+1 \over 2} .
\ea
\eeq

Let us look now at the convergence guarantees for method (\ref{met-DEGM}) with tensor step (\ref{def-BTStep}). Defining $V(x) = \nabla f(x)$, we get
\beq\label{eq-DMain}
\ba{rcl}
\sum\limits_{i=1}^t a_i \psi(x_i) + B_t & \refLE{eq-DEM} & \Psi_t^* = \min\limits_{x \in \dom \psi} \Big[ \half \| x - x_0 \|^2 + \sum\limits_{i=1}^t a_i[\la \nabla f(x_i), x - x_i \ra] + A_t \psi(x) \Big]\\
\\
& \leq & \half \| x^* - x_0 \|^2 + \sum\limits_{i=1}^t a_i[\la \nabla f(x_i), x^* - x_i \ra] + A_t \psi(x^*)\\
\\
& \leq & \half \| x^* - x_0 \|^2 + A_t F(x^*) -  \sum\limits_{i=1}^t a_i f(x_i).
\ea
\eeq
Hence, $B_t + A_t[\tilde F_t - F_*] \leq \half \| x^* - x_0 \|^2$,
and we conclude that for this method, the same rates of convergence (\ref{eq-RateGP}), (\ref{eq-RateFP}), and (\ref{eq-RateFP*}) are valid.

Its accuracy certificate need to be defined in a different way, namely
\beq\label{def-DCertF}
\ba{rcl}
\Delta^{\psi}_t & = & {1 \over A_t} \max\limits_{x \in \dom \psi} \Big\{ \sum\limits_{i=1}^t a_i [ \la \nabla f(x_i), x_i - x \ra + \psi(x_i) - \psi(x)]: \; \| x - x_0 \| \leq R_0 \Big\}\\
\\
& \geq & {1 \over A_t} \max\limits_{x \in \dom \psi} \Big\{ \sum\limits_{i=1}^t a_i [ F(x_i) - F(x)]: \; \| x - x_0 \| \leq R_0 \Big\} \; = \; \tilde F_t - F_*.
\ea
\eeq
On the other hand,
$$
\ba{rcl}
\sum\limits_{i=1}^t a_i \psi(x_i) + B_t & \refLE{eq-DMain} &  \half R_0^2 + \min\limits_{x \in \dom \psi} \Big\{ \sum\limits_{i=1}^t a_i[\la \nabla f(x_i), x - x_i \ra + \psi(x)]:\; \| x - x_0 \| \leq R_0 \Big\}\\
\\
& = & \half R_0^2 + \sum\limits_{i=1}^t a_i \psi(x_i) - A_t \Delta^{\psi}_t.
\ea
$$
Thus, $B_t + A_t \Delta^{\psi}_t \leq \half R_0^2$. Hence, this certificate satisfies the bounds (\ref{eq-RateCF}) and (\ref{eq-RateCF*}).

Finally, for method (\ref{met-PREGM}) with step (\ref{def-BTStep}), the only performance guarantee is as follows:
$$
\ba{rcl}
B_t & \leq & \half \| x_0 - x^* \|^2, \quad t \geq 1.
\ea
$$
Hence,we can justify only the bounds (\ref{eq-RateGP}) and (\ref{eq-RateGP*}) for the norm of the gradient.

To conclude this section, note that the rate of convergence of method (\ref{met-EPrimal}), (\ref{def-BTStep}) for the norms of reduced gradients can be improved by a simple switching strategy. 

Let us fix a number of steps $t \geq 1$. The total number of steps in our method will be $N = 2t$. It consists of two stages.
\beq\label{eq-EPrim}
\ba{rl}
\mbox{a)} & \mbox{Compute points $\{x_i\}_{i=1}^t$ and coefficients $\{ a_i \}_{i=1}^t$ by method (\ref{met-EPrimal}), (\ref{def-BTStep}).}\\
\mbox{b)} & \mbox{Define $y_0 = {1 \over A_t} \sum\limits_{i=1}^t a_i x_i$. For $i=0, \dots, t-1$, iterate $y_{i+1} = T_M^p(y_i)$.}
\ea
\eeq
Denote $G_i = \| F'(y_i \|_*$, $1 \leq i \leq t$, and $G_N^* = \min\limits_{1 \leq i \leq t} G_i$. 

Let us choose $M = M^*_p$. Note that
$$
\ba{rcl}
{\| x_0 - x^* \|^{p+1} \over p+1} \left({1 \over \gamma^*_p}\right)^p \left(p-1 \over p+1\right)^{p-1 \over 2} \left(1 \over t \right)^{p+1 \over 2} & \refGE{eq-RateFP} & F(y_0) - F_*.
\ea
$$
On the other hand, for all $i=0, \dots, t-1$, we have
$$
\ba{rcl}
F(y_i) - F(y_{i+1}) & \geq & \la F'(y_{i+1}), y_i - y_{i+1} \ra \; \refGE{eq-VIStep} \; \gamma_p^* G_{i+1}^{p+1 \over p}.
\ea
$$
Hence, $F(y_0) - F_* \geq t \cdot \gamma_p^* \left(G^*_N \right)^{p+1 \over p}$, and we get
\beq\label{eq-GFast}
\ba{rcl}
G_N^* & \leq & \left( {\| x_0 - x^* \| \over \gamma_p^*} \right)^p \left({1 \over p+1} \right)^{p \over p+1} \left( p-1 \over p+1 \right)^{p(p-1) \over 2(p+1)} \left( {2 \over N} \right)^{p(p+3) \over 2(p+1)} \\
\\
& \leq & {2 L_p \| x_0 - x^* \|^p \over p!} \left( {2 \over N} \right)^{p(p+3) \over 2(p+1)}.
\ea
\eeq

\section{High-order methods for variational inequalities}\label{sc-Var}
\SetEQ

Consider now the general problem of Variational Inequality (\ref{prob-CVar}). 
We will measure quality of the approximate solution $\bar x \in \dom \psi$ by the following {\em composite merit function}:
$$
\ba{rcl}
\mu_{\psi}(\bar x) & = & \psi(\bar x) + \max\limits_{x \in \dom  \psi} \Big[ \la V(x), \bar x - x \ra - \psi(x) \Big],
\ea
$$
Clearly, for all $\bar x \in \dom \psi$, we have $\mu_{\psi}(\bar x) \geq 0$, and $\mu_{\psi}(x^*) = 0$ if and only if $x^* \in \dom \psi$ is a solution to (\ref{prob-CVar}). 

We are going to measure distances in the Euclidean norm $\| \cdot \| = \| \cdot \|_B$, assuming that the set $\dom \psi$ is bounded:
$$
\ba{rcl}
\| x - x_0 \| & \leq & R_0, \quad \forall x \in \dom \psi.
\ea
$$
Thus, for any feasible point $x$, we have
$\beta(x_0,x) = \half \| x - x_0 \|^2 \leq \half R_0^2$.

Let us estimate the quality of the following output of our methods (\ref{met-EPrimal}) and (\ref{met-DEGM}):
\beq\label{def-Out}
\ba{rcl}
\bar x_t & = & {1 \over A_t} \sum\limits_{i=1}^t a_i x_i.
\ea
\eeq
We can do this, using the following {\em computable} accuracy certificate
\beq\label{def-Cert}
\ba{rcl}
\Delta^{V}_t & \Def & {1 \over A_t} \max\limits_{x \in \dom \psi} \sum\limits_{i=1}^t a_i [\la V(x_i), x_i - x \ra + \psi(x_i) -\psi(x)]\\
\\
& \refGE{eq-Mon} & \psi(\bar x_t) + {1 \over A_t} \max\limits_{x \in \dom \psi} \sum\limits_{i=1}^t a_i [\la V(x), x_i - x \ra -\psi(x)]\\
\\
& \refEQ{def-Out} & \psi(\bar x_t) + \max\limits_{x \in \dom \psi} \Big[ \la V(x), \bar x_t - x \ra -\psi(x) \Big] \; = \; \mu_{\psi}(\bar x_t).
\ea
\eeq

For the primal method (\ref{met-EPrimal}), for any $x \in \dom \psi$, we have
$$
\ba{rcl}
\half \| x - x_0 \|^2 & \stackrel{(\ref{eq-EPrimal}),(\ref{eq-RGLow})}{\geq} & \half \| x - v_t \|^2 + B_t + \sum\limits_{i=1}^t a_i [ \la V(x_i), x_i - x \ra + \psi(x_i) - \psi(x)].
\ea
$$
Hence, taking $x=x^*$, we get
\beq\label{eq-BoundVB}
\ba{rcl}
\half \| v_t - x^*  \|^2 + B_t & \refLE{prob-CVar} & \half \| x_0 - x^* \|^2.
\ea
\eeq
On the other hand, by the uniform bound on the size of variables, we have
\beq\label{eq-RateVP}
\ba{rcl}
B_t + A_t \Delta^{V}_t & \leq & \half R_0^2. 
\ea
\eeq

For the dual method (\ref{met-DEGM}), we have
\beq\label{eq-TDual}
\ba{rcl}
\sum\limits_{i=1}^t a_i \psi(x_i) + B_t \; \refLE{eq-DEM} \; \Psi^*_t 
& = & \min\limits_{x \in \dom \psi} \Big\{ \half \| x - x_0 \|^2 + \sum\limits_{i=1}^t a_i [ \la V(x_{i}), x - x_i \ra + \psi(x) ] \Big\}.
\ea
\eeq
Since function $\Psi_t(\cdot)$ in the dual method is strongly convex with convexity parameter one, this means that
$$
\ba{rcl}
\sum\limits_{i=1}^t a_i \psi(x_i) + B_t + \half \| x^* - v_t \|^2 & \leq & \half \| x^* - x_0 \|^2 + \sum\limits_{i=1}^t a_i [ \la V(x_{i}), x^* - x_i \ra + \psi(x^*) ] \\
\\
& \refLE{prob-CVar} & \half \| x^* - x_0 \|^2 + \sum\limits_{i=1}^t a_i \psi(x_i),
\ea
$$
and this is the relation (\ref{eq-BoundVB}). On the other hand, we have 
$$
\ba{rcl}
\sum\limits_{i=1}^t a_i \psi(x_i) + B_t
& \refLE{eq-TDual} & \half R_0^2 - A_t \Delta^{V}_t + \sum\limits_{i=1}^t a_i \psi(x_i).
\ea
$$
Thus, inequality (\ref{eq-RateVP}) is valid for the dual method too.
It remains to justify a certain rate of growth of the coefficients $A_t$.

Let us fix some $p \geq 0$ and assume that $(p+1)$st derivative of $V(\cdot)$ is uniformly bounded:
\beq\label{ed-DSBound}
\ba{rcl}
\| D^{p+1}V(x) \|^0 & \leq & \hat M_{p+1}, \quad x \in \dom \psi.
\ea
\eeq
Consider the following operator
$$
\ba{rcl}
G_{x,p}^M(y) & = & {\cal T}_{x,p}^V(y) + M \nabla d_{p+2}(y-x), \quad y \in \dom \psi.
\ea
$$
In our notation, the case $p=0$ corresponds to the {\em zero-order methods}, which use no derivatives of the operator $V(\cdot)$. Then the operator ${\cal T}_{x,0}^V(\cdot) \equiv V(x) \in \E^*$ is constant.

\BT\label{th-Pos}
For any $M \geq { 1 \over p!} \hat M_{p+1}$, operator $G_{x,p}^M(\cdot)$ is monotone on $\dom \psi$.
\ET
\proof
Since $D {\cal T}_{x,p}^V(y) \refEQ{eq-V2Div} {\cal J}_{x,p}^V(y)$, for any $h \in \E$, we have
$$
\ba{rcl}
\la D G_{x,p}^M(y)h,h \ra & = & \la {\cal J}_{x,p}^V(y)h,h \ra + M \la \nabla^2 d_{p+2}(y-x) h,h \ra \\
\\
& \refGE{eq-DPChar} & M \| x - y \|^{p} \| h \|^2 + \la \Big[{\cal J}_{x,p}^V(y) - DV(y)\Big]h,h \ra + \la DV(y)h,h \ra\\
\\
& \refGE{eq-V2Div} & M \| x - y \|^{p} \| h \|^2 - {\hat M_{p+1} \over p!} \| y - x \|^p \| h \|^2 \geq 0. \QR
\ea
$$

Let us analyze now the Basic Tensor Step for variational inequalities.
It consists in solving the following auxiliary problem:
\beq\label{eq-CondVI}
\fbox{$\ba{c}\\
\mbox{Find $x_+ \Def V_M^p(v) \in \dom \psi$ satisfying the following condition:}\\
\\
\la G_{v,p}^M(x_+), y - x_+ \ra + \psi(y) \geq \psi(x_+), \; \forall y \in \dom \psi.\\
\\
\ea$}
\eeq
Thus, for computing $x_+$, we need to solve a variational inequality with monotone operator. However, for that it is not necessary to call oracle of the initial operator $V(\cdot)$.

In view of Lemma \ref{lm-Cut2}, for point of the form $x_+ = V_M^p(v)$ with $v \in \dom \psi$, we can define 
\beq\label{eq-VPVar}
\ba{rcl}
V_{\psi}(x_+) & = & V(x_+) - G_{v,p}^M(x_+).
\ea
\eeq
Note that $G_{x^*,p}(x^*) =V(x^*)$. Therefore, $x^* \refEQ{prob-CVar1} V^p_M(x^*)$, and we conclude that
\beq\label{eq-VZero}
\ba{rcl}
V_{\psi}(x^*) & = & 0.
\ea
\eeq

Let us show how the Reduced Gradient looks like in some simple situation.
\BE\label{ex-VPsi}
Let $Q \subseteq \E$ be a closed convex set, $\psi(x) = \Ind_Q(x)$, and $p=0$. Then 
$$
\ba{rcl}
G_{v,0}^M(y)& = & V(v) + M B(y - v), \quad y_+ \; \Def \; v - {1 \over M} B^{-1} V(v), \quad x_+ \; = \; V^0_M(v) \; = \; \pi_Q(y_+).
\ea
$$
Thus, we get 
$$
\ba{rcl}
V_{\psi}(x_+) & \refEQ{eq-VPVar} & V(x_+) - V(v) - M B (x_+ - v) \; = \; V(x_+) + M d(x_+), \\
\\
d(x_+) & \Def & B \left(y_+ - x_+ \right). 
\ea
$$
If $Q=\E$, then $V_{\psi}(x_+) = V(x_+)$. However, in the constrained case, our cutting plane differs from the classical extragradient \cite{Korp, Arik}. If $y_+ \not\in Q$, then the value $V(x_+)$ is corrected by an additive term, which belongs to the normal cone at point $x_+ \in \partial Q$. This correction ensures vanishing value of $V_{\psi}(x_+)$ as $x_+ \to x^*$, explaining our choice of the name "reduced gradient". It corresponds to the common practice of using such a terminology in Optimization \cite{Wolfe}.
\qed
\EE

Let us find now a lower bound for the universal stepsize (\ref{def-ASize}). Note that
$$
\ba{rcl}
\la V_{\psi}(x_+), v - x_+ \ra & = & \la V(x_+) - {\cal T}_{v,p}^V(x_+) + M \nabla d_{p+2} (v - x_+), v - x_+ \ra \\
\\
& \refGE{eq-VDiv} & M \| v - x_+ \|^{p+2} - {\hat M_{p+1} \over (p+1)!} \| v - x_+ \|^{p+2} \\
\\
& = & \left(M - {\hat M_{p+1} \over (p+1)!} \right) \| v - x_+ \|^{p+2}.
\ea 
$$
On the other hand
$$
\ba{rcl}
\| V_{\psi}(x_+) \|_* & = & \| V(x_+) - {\cal T}_{v,p}^V(x_+) + M \nabla d_{p+2} (v - x_+) \|_*\\
\\
& \refLE{eq-VDiv} & \left(M + {\hat M_{p+1} \over (p+1)!} \right) \| v - x_+ \|^{p+1}.
\ea
$$
Thus, we get
\beq\label{eq-VSizeL}
\ba{rcl}
\la V_{\psi}(x_+), v - x_+ \ra & \geq & \left(M - {\hat M_{p+1} \over (p+1)!} \right) \Big[ \left(M + {\hat M_{p+1} \over (p+1)!} \right)^{-1} \| V_{\psi}(x_+) \|_* \Big]^{p+2 \over p+1}\\
\\
& \Def & \hat \gamma_p  \| V_{\psi}(x_+) \|_*^{p+2 \over p+1}.
\ea
\eeq

This means that we can estimate the coefficients in methods (\ref{met-EPrimal}) and (\ref{met-DEGM}) as follows:
\beq\label{eq-LCoef}
\ba{rcl}
a_{t+1} & = & {\la V_{\psi}(x_{t+1}), v_t - x_{t+1} \ra \over \| V_{\psi}(x_{t+1}) \|_*^2} \; \geq \; \hat \gamma_p \| V_{\psi}(x_{t+1}) \|_*^{-{p \over p+1}},\\
\\
b_{t+1} & = & {\la V_{\psi}(x_{t+1}), v_t - x_{t+1} \ra^2 \over 2 \| V_{\psi}(x_{t+1}) \|_*^2} \; \geq \; \half \hat \gamma^2_p \| V_{\psi}(x_{t+1}) \|_*^{{2 \over p+1}}.
\ea
\eeq

For $i \geq 1$, denote $g_i = \| V_{\psi}(x_{i}) \|_*$, and $g^*_t = \min\limits_{1 \leq i \leq t} g_i$. Then from inequality (\ref{eq-BoundVB}), we have
\beq\label{eq-VRateG}
\ba{rcl}
g^*_t & \leq & \left({1 \over t}\right)^{p+1 \over 2} \left[ {1 \over \hat \gamma_p} \| x_0 - x^* \| \right]^{p+1}, \quad t \geq 1,
\ea
\eeq
where 
\beq\label{def-GP}
\ba{rcl}
\hat \gamma_p & = & \left(M - {\hat M_{p+1} \over (p+1)!} \right)  \left(M + {\hat M_{p+1} \over (p+1)!} \right)^{-{p+2 \over p+1}}. 
\ea
\eeq
The minimal value of the constant $\hat \gamma_p = \hat \gamma_p(M)$ is attained for $M = \hat M^*_{p+1} = {2p+3 \over (p+1)!} \hat M_{p+1}$. That is
$$
\ba{rcl}
\hat \gamma^*_p & = & {p+1 \over p+2} \left[ (p+1)! \over 2(p+2) \hat M_{p+1} \right]^{1 \over p+1}.
\ea
$$
In this case, since $\left(1 + {1 \over p+1} \right)^{p+1} \leq e$, we get the following bound:
\beq\label{eq-VRateG*}
\ba{rcl}
g^*_t & \leq & {2e(p+2) \over (p+1)!}\left({1 \over t}\right)^{p+1 \over 2} \hat M_{p+1} \| x_0 - x^* \|^{p+1}, \quad t \geq 1.
\ea
\eeq

Let us obtain now the rate of convergence in terms of the accuracy certificate. In view of inequality (\ref{eq-RateVP}), we have
$$
\ba{rcl}
\half R_0^2 & \refGE{eq-LCoef} & \hat \gamma_p \sum\limits_{i=1}^t \Big[ \half \hat \gamma_p g_i^{2 \over p+1} +  g_i^{-{p \over p+1}} \Delta^{V}_t \Big] \; \refGE{eq-Tech} \;
\hat \gamma_p \; t \cdot {p+2 \over 2} \left(\hat \gamma_p \over p\right)^{p \over p+2}  (\Delta^{V}_t)^{2 \over p+2}.
\ea
$$

Thus, we have proved the following theorem.
\BT\label{th-RateVI}
Let $p \geq 0$ and $\hat M_{p+1} < + \infty$. Let us choose $M \geq {1 \over p!}\hat M_{p+1}$. Then both methods (\ref{met-EPrimal})and (\ref{met-DEGM}) with the essential step (\ref{eq-CondVI}) and the cutting plane defined by (\ref{eq-VPVar}) ensure the following rates of convergence:
\beq\label{eq-RateVI}
\ba{rcl}
g^*_t & \leq & \left({1 \over t}\right)^{p+1 \over 2} \left[ {1 \over \hat \gamma_p} \| x_0 - x^* \| \right]^{p+1}, \\
\\
\mu_{\psi}(\bar x_t) \; \leq \; \Delta^{V}_t  & \leq &  \left(1 \over t \right)^{p+2 \over 2} \left(1 \over \hat \gamma_p \right)^{p+1} {R_0^{p+2} \over p+2}, \quad t \geq 1,
\ea
\eeq
where $\hat \gamma_p$ is defined by (\ref{def-GP}).
\ET

For the optimal value of $M = \hat M^*_{p+1} = {2 p + 3 \over p!} \hat M_{p+1}$, these bounds have the following form:
\beq\label{eq-RateVI*}
\ba{rcl}
g^*_t & \leq & {2e(p+2) \over (p+1)!}\left({1 \over t}\right)^{p+1 \over 2} \hat M_{p+1} \| x_0 - x^* \|^{p+1},\\
\\
\Delta^V_t & \leq & {2e \over (p+1)!}\left({1 \over t}\right)^{p+2 \over 2} \hat M_{p+1} R_0^{p+2}, \quad p \geq 0.
\ea
\eeq

Recall that the latter performance guarantee is unimprovable \cite{ABJS}. To the best of our knowledge, the lower bounds for the norms of reduced gradients are not known yet.

\section{Methods for uniformly monotone CVI}\label{sc-Uni}
\SetEQ

Let us assume now that the operator $V(\cdot)$ in problem (\ref{prob-CVar}) is uniformly monotone:
\beq\label{def-UMon}
\ba{rcl}
\la V(x) - V(y), x - y \ra & \geq & \sigma_{p+2} \| x - y \|^{p+2}, \quad x, y \in \dom \psi,
\ea
\eeq
with the degree of monotonicity $p+2$, $p \geq 0$, and parameter $\sigma_{p+2} > 0$. The important case $p=0$ corresponds to {\em strong monotonicity}. For simplicity, we assume that the norm is Euclidean: 
$$
\ba{rcl}
\beta(x_0,x) & = & \half \| x - x_0 \|^2_B.
\ea
$$

Let us look at some variant of the primal method (\ref{met-EPrimal}) with the essential step (\ref{eq-CondVI}). 
Condition (\ref{def-UMon}) can be used for deriving its linear rate of convergence for the distance to the optimal solution $x^*$ (which is unique now). Let us introduce the condition number $\kappa_{p+1} \Def {\sigma_{p+2} \over \hat M_{p+1}}$. 

\beq\label{met-Prop}
\ba{|l|}
\hline \\
\mbox{{\bf 0.} Choose $x_0 \in \dom \psi$ and $p \geq 0$. Set $v_0 = x_0$.}\\
\\
\mbox{\bf 1. $t$-th iteration ($t \geq 0$).}\\
\\
\mbox{{\bf a).} For $M \geq {1 \over p!} \hat M_{p+1}$, compute $x_{t+1} = V^p_{M}(v_t)$ and define $V_{\psi}(x_{t+1})$ by (\ref{eq-VPVar}).}\\
\\
\mbox{{\bf b).} Set $a_{t+1} = a(v_t,x_{t+1})$ and compute $\hat v_{t+1} = \prox_{v_t}^{a_{t+1}}(V_{\psi}(x_{t+1}))$.}\\
\\
\mbox{{\bf c).} For 
$\alpha_p = (p+2) \hat \gamma_p \left({1 \over p} \hat \gamma_p \right)^{p \over p+2} \sigma_{p+2}^{{2 \over p+2}}$, define $v_{t+1} = {1 \over 1 + \alpha_p} (\hat v_{t+1} + \alpha_p x_{t+1})$.}\\
\\
\hline
\ea
\eeq
Note that for $\alpha_p = 0$, method (\ref{met-Prop}) coincides with method (\ref{met-EPrimal}), (\ref{eq-CondVI}). 

\BT\label{th-Prop}
Let $\kappa_{p+1} > 0$ for some $p \geq 0$. Then for the sequence of points $\{ x_t \}_{t \geq 0}$, generated by method (\ref{met-Prop}), we have
\beq\label{eq-PRate}
\ba{rcl}
\| v_t - x^* \|_B & \leq &  \left( 1 \over 1 + \alpha_p \right)^{t/2} \| x_0 - x^* \|_B , \quad t \geq 0.
\ea
\eeq
\ET
\proof
In the proof of Theorem \ref{th-EPrimal}, let us take $x=x^*$. Then for any $t \geq 0$, we have
$$
\ba{c}
\half \| v_t - x^* \|^2_B - \half \| \hat v_{t+1} - x^* \|^2_B \; \geq \; b_{t+1} + a_{t+1} V_{t+1}(x^*) \;
\refEQ{def-Obj} \; b_{t+1} + a_{t+1} \la V_{\psi}(x_{t+1}), x_{t+1} - x^* \ra \\
\\
\refGE{eq-VPsi} \; b_{t+1} + a_{t+1} \la V_{\psi}(x^*), x_{t+1} - x^* \ra + a_{t+1} \la V(x_{t+1}) - V(x^*), x_{t+1} - x^* \ra\\
\\
\stackrel{(\ref{eq-VZero}), (\ref{def-UMon})}{\geq} \; b_{t+1} + a_{t+1} \sigma_{p+2} \| x_{t+1} - x^* \|^{p+2}.
\ea
$$
Denoting $g_{t+1} = \| V_{\psi}(x_{t+1}) \|_*$ and $r_{t+1} = \| x_{t+1} - x^* \|_B$, and using the lower bounds (\ref{eq-LCoef}), we can apply Lemma \ref{lm-Tech} with $\sigma = p+1$ and estimate the right-hand side of this inequality as follows:
$$
\ba{rcl}
\hat \gamma_p  \left( \half \hat \gamma_p g_{t+1}^{2 \over p+1} + \sigma_{p+2} g_{t+1}^{- {p \over p+1}} r_{t+1}^{p+2} \right) & \geq & {p+2 \over 2} \hat \gamma_p \left({1 \over p} \hat \gamma_p \right)^{p \over p+2} \sigma_{p+2}^{{2 \over p+2}}r_{t+1}^2 \; \Def \; \half \alpha_p r_{t+1}^2.
\ea
$$
Hence, $\| v_t - x^* \|^2_B \geq  \| \hat v_{t+1} - x^* \|^2_B + \alpha_p \| x_{t+1} - x^* \|^2_B \geq (1+\alpha_p) \| v_{t+1} - x^* \|^2_B$.
\qed

For the optimal choice $M = \hat M^*_{p+1} = {2 p + 3 \over p!} \hat M_{p+1}$, we have $\alpha_p = O\left(\kappa_p^{2 \over p+2}\right)$.
Hence, in order to ensure $\| v_t - x^* \|_B \leq \epsilon$, method (\ref{met-Prop}) needs
$O \left( \left[ {\hat M_{p+1} \over \sigma_{p+2}} \right]^{2 \over 2+p} \ln {\| x_0 - x^* \| \over \epsilon} \right)$
calls of oracle.

\section{Conclusion}\label{sc-Conc}
\SetEQ

This paper can be seen as an alternative to the ideas presented in two recent paper \cite{ABJS} and \cite{LJ}, where the authors justified two different high-order extragradient methods for solving monotone Variational Inequalities. In our terminology, method in \cite{ABJS} is an analogue of the Primal Extragradient Method (\ref{met-EPrimal}). Method in \cite{LJ}, which can be seen as an extension of the Dual Extrapolation Method \cite{DExtra}, is similar to the Dual Extragradient Method (\ref{met-DEGM}). For Variational Inequalities, all three methods give similar performance guarantees \footnote{For being consistent with description of the minimization schemes, our definition of order of the methods for VIs is shifted by one as compared to \cite{ABJS,LJ}. Thus, for us the classical Extragradient Method \cite{Korp, Arik} is the {\em zeroth-order scheme} since it does not use derivatives of the operator.}.

Let us discuss the aspects in which our approach is different from that of our predecessors.
\begin{itemize}
\item
Our methods are based on the concept of Reduced Gradient, which is different from the classical Extra-Gradient approach (see discussion in Example \ref{ex-VPsi}).
\item
In both papers \cite{ABJS,LJ}, the size of the essential step is related to some power of the norm of the extragradient step, which is prescribed by the order of Taylor polynomial in use. Hence, their results assume an involvement of Taylor polynomial of certain degree.

Our approach is more flexible since the big size of the Universal Step (\ref{def-ASize}) can  be achieved by other strategies for computing the cutting plane (\ref{eq-Cut}).
\item
Our method for uniformly monotone operators is different from the previously known schemes (see \cite{LJ}). It does not need restarts and does not require preliminary knowledge of the size of the feasible set.
\item
All our methods benefit from the hot-start capabilities, which can be useful even for optimization schemes.
\item
We obtain also the rates of convergence for the norms of reduced gradients, which seems to be new measures of optimality of the approximate solutions to CVIs.
\item
Our results are related to Composite Variational Inequality Problem (\ref{prob-CVar}), which seems to be a new and useful problem setting. 
\end{itemize}

\end{document}